\documentclass[10pt]{amsart}

\usepackage{amsmath,mathtools,amsthm,amssymb,dsfont,mathrsfs,color,graphicx,amssymb,bm}
\usepackage{graphicx}
\usepackage{placeins}
\usepackage{algorithm,algpseudocode}
\usepackage{caption}
\usepackage{esint}
\usepackage[svgnames]{xcolor} 
\usepackage[colorlinks,citecolor=red,pagebackref,hypertexnames=false,breaklinks]{hyperref}
\usepackage{pdfsync}

\newtheorem{theorem}{Theorem}

\newcommand{\supp}[1]{{\rm supp}\,\,#1}
\newcommand{\singsupp}[1]{{\rm sing}\,\,{\rm supp}\,\,#1}
\newcommand{\N}{\mathcal{N}}
\newcommand{\norm}[1]{\left\|#1\right\|}

\newcommand{\bea}{\begin{eqnarray*}}
\newcommand{\eea}{\end{eqnarray*}}
\newcommand{\bean}{\begin{eqnarray}}
\newcommand{\eean}{\end{eqnarray}}
\newcommand{\be}{\begin{equation}}
\newcommand{\ee}{\end{equation}}
\newcommand{\ba}{\begin{array}}
\newcommand{\ea}{\end{array}}

\title[Localization of moving sources]{Localization of moving sources:\\uniqueness, stability, and Bayesian inference}

\author{S\'ara Wang}
\thanks{Department of Applied Mathematics and Computer Science, Technical University of Denmark, Danemark. Email: s144275@student.dtu.dk}
\author{Mirza Karamehmedovi\'c}
\thanks{Department of Applied Mathematics and Computer Science, Technical University of Denmark, Denmark. Email: mika@dtu.dk. M. K. was supported by The Villum Foundation (grant no. 25893).}
\author{Faouzi Triki}
\thanks{Laboratoire Jean Kuntzmann, UMR CNRS 5224, Universit\'e Grenoble-Alpes, 700 Avenue Centrale, 38401 Saint-Martin-d'H\`eres, France. E-mail: faouzi.triki@univ-grenoble-alpes.fr. F. T. was supported by the grant ANR-17-CE40-0029 of the French National Research Agency ANR (project MultiOnde)}

\date{\today}

\makeatletter
\newcommand\ackname{Acknowledgements}
\if@titlepage
  \newenvironment{acknowledgements}{%
      \titlepage
      \null\vfil
      \@beginparpenalty\@lowpenalty
      \begin{center}%
        \bfseries \ackname
        \@endparpenalty\@M
      \end{center}}%
     {\par\vfil\null\endtitlepage}
\else
  
\fi
\makeatother

\begin{document}

\maketitle

\begin{abstract}
We consider the subsonic moving point source problem for the scalar wave equation in $\pmb{R}^3$, proving a regularity result for the direct problem,  and uniqueness and stability results for the inverse problem. We then present and investigate numerically a Bayesian framework for the inference of the source trajectory and intensity from wave field measurements. The framework employs Gaussian process priors, the pre-conditioned Crank-Nicholson scheme with Markov Chain Monte Carlo sampling, and conditioning on functionals to include prior information on the source trajectory.
\end{abstract}


\section{Introduction}\label{sec:introduction}

This work concerns the moving source problem for the scalar wave equation in $\pmb{R}^3$,\begin{equation}\label{eqn:system}
\left\{
\begin{array}{rcll}
c^{-2}\partial^2_tu-\Delta u&=&q(t)\delta_{p(t)}\in C^m(]0,T[,\mathcal{E}'^0(\pmb{R}^3)),&\\
u(0+,x)&=&0,\quad x\in\pmb{R}^3,\\
\partial_tu(0+,x)&=&0,\quad x\in\pmb{R}^3.
\end{array}\right.
\end{equation}
Here $m\in\pmb{N}$ is fixed, $c\in\pmb{R}_+$ is the constant wave speed, $T\in\pmb{R}_+$ is a fixed final time, $p\in C^{m+1}(]0,T[, \Omega_0)$ is the source trajectory (with $\Omega_0\subset\pmb{R}^3$ an open, bounded, simply connected domain), $q\in C^m(]0,T[)$ is the nonnegative-valued source intensity,  $\delta_a$ is the Dirac delta supported at $a\in\pmb{R}^3$,
and $\mathcal{E}'^0(\pmb{R}^3)$ is the space of compactly supported distributions. We consider subsonic sources, that is, we assume
\bean \label{subsonic}
|p'(t)|<c,\quad t\in]0,T[.
\eean 
The fundamental solution $E_+\in C^{\infty}(]0,T[,\mathcal{E}'^0(\pmb{R}^3))$ of the wave operator, given by~\cite[Section 6.2, Eqs. (6.2.4)' and (6.2.6)']{HI}
\[
E_+(t)(\phi)=\frac{c^2t}{4\pi}\int_{\omega\in S^2}\phi(ct\omega)d\omega
,\quad t\in]0,T[,\,\,\phi\in C^{\infty}(\pmb{R}^3),
\]
satisfies
\begin{equation*}
\left\{
\begin{array}{rcl}
(c^{-2}\partial^2_t-\Delta)E_+(t)&=&\delta_0,\quad t\in]0,T[,\\
E_+(0+)&=&0,\\
\partial_tE_+(0+)&=&c^2\delta_0,
\end{array}\right.
\end{equation*}
and the solution~\cite[Theorem 6.2.4]{HI}
\begin{align}\label{eqn:solution_integral}
u(t,\cdot)&=\int_0^tE_+(t-s)\ast(q(s)\delta_{p(s)})ds\nonumber\\&=\frac{c^2}{4\pi}\int_0^t q(s)(t-s)\int_{\omega\in S^2}\delta_{c(t-s)\omega +p(s)}d\omega ds\in C^{\infty}(]0,T[,\mathcal{E}'^0(\pmb{R}^3)),
\end{align}
of~\eqref{eqn:system} is the field radiated by the point source.
It is well-known~\cite{Jackson} that $u(t,x)$ can be expressed using the Li\'enard-Wiechert retarded potential
\begin{equation}\label{eqn:solution}
u(t,x)=\begin{cases}(c/4\pi)q(\tau(x,t))|x-p(\tau(x,t))|^{-1}h(x,\tau(x,t))^{-1},&|x-p(0)|<ct,\\0&{\rm otherwise,}\end{cases}
\end{equation}
where $\tau(x,t)$ is given implicitly by 
\bean \label{tau}
c(t-\tau(x,t))=|x-p(\tau(x,t))|, 
\eean
and
\bean \label{h}
h(x,\tau)=c-\frac{x-p(\tau)}{|x-p(\tau)|}p'(\tau)\quad\text{for}\,\,x\neq p(\tau).
\eean
Notice that since $p$  satisfies the inequality \eqref{subsonic}, we have
\[
\frac{d}{d\tau}\left(\tau\mapsto t-c^{-1}|x-p(\tau)|\right)\le c^{-1}|p'(\tau)|<1,
\]
so the fixed point problem \eqref{tau} has a unique solution.

Now let $\Omega$ be a bounded and connected domain with a $C^2$ boundary $\partial \Omega$, and such that $\overline{\Omega}_0\subset \Omega$. We denote by $\nu(x)$ the unit outward normal vector at $x\in \partial \Omega$. In the following we measure the trace of the solution $u$ at the {\it measurement surface} $\Gamma\subseteq\partial\Omega$ that is open in the induced topology on $\partial\Omega$, and that is not included in any plane intersecting $\Omega$. The inverse moving source problem we consider is to reconstruct $p$ and $q$, given a measurement $u|_{[0,T]\times\Gamma}$; the associated forward problem is to find $u$ given $p$ and $q$. We must {\it turn off} the source at a time $T_0<T$ (i.e., $q(t)=0$ for $t\in[T_0,T]$) to have any chance of uniqueness of solution of the inverse problem, due to the finite speed of wave propagation. We consider $T$ large enough to allow the information about $p$ and $q$ for $t\in [0, T_0]$ to propagate towards the measurement set $\Gamma$, that is,  
\bean \label{timeassumption}
T>T_0+\sup\{|x-y|/c,\,\,x\in\Gamma,\,\,y\in\Omega_0\} := T^*.
\eean
We further assume that
the set $\{t\in[0,T_0],\,\,q(t)=0\}$ has zero measure. This prevents the situation where two trajectories are indistinguishable because the source does not radiate where the trajectories are not aligned. Finally, if $\Gamma$ is included in a plane $P$ intersecting $\Omega$ then one is unable to distinguish the field measurements from pairs of sources constructed using the same intensity, and by translating a trajectory confined in $P\cap\Omega$ to two different planes parallel to $P$ and at the same orthogonal distance from $P$. (More generally, trajectories mirrored across $P$ and having the same amplitude profile $q$ would also produce the same field measurements.)

The inverse moving source problem for the wave equation  was treated in many recent works~\cite{ohe,nakaguchi,hu,hu2, AlJebawy}. Different approaches have been considered for solving it. In \cite{el2001determination, ohe2011real}, the authors have applied  direct algebraic methods for reconstructing stationary point sources from  boundary measurement. Later, in \cite{ohe}, this  approach was extended to the problem of moving sources with a strong assumption on their trajectories. A different algebraic method  for the identification of a single moving source using measurements of the retarded potential and all its derivatives at a single observation point was provided in  \cite{nakaguchi2012algebraic}. Optimization techniques were also applied in solving the inverse  problem (see for instance \cite{bruckner2000determination,komornik2002upper,komornik2005estimation,rashedi2015stable}). Recently a Lipschitz stability estimate was derived for recovering a single moving point source from the knowledge of the field at six well-chosen points on the boundary assuming that the intensity  is a known constant  \cite{AlJebawy}. Notice that when the point source is stationary  the inverse problem is well studied, and H\"older type stability estimates are derived using a single boundary measurement  \cite{baorecovering,bao2021recovering}.  \\

Our main theoretical results concern the regularity of solution of the forward problem, and the uniqueness and stability of solution of the inverse problem. 

\begin{theorem}\label{thm:regularity} Let $u$ be the solution  of the system~\eqref{eqn:system}. Then
 $$u\in C^{\infty}(]0,T[,H^{m-1/2-\varepsilon}(\pmb{R}^3)),$$ for every strictly positive $\varepsilon$. Moreover, \[u\in C^{\infty}(]0,T[,C^{\infty}(\pmb{R}^3\setminus\cup_{s=0}^{T_0}\overline{B_{p(s)}(c(T_0-s))})),\] and $u(t,\cdot)$ is $C^{\infty}$ in \[\Bigl\{|x-p(s)|<c(t-s),\,\,0\le s\le\min\{t,T_0\}\Bigr\}\setminus p([0,T_0]).\]
 Finally, 
 \bea
 u(t, \cdot)|_\Gamma = 0 \quad \textrm{for  } t \in ]T^*, +\infty[,
 \eea
 with $T^* = T_0+\sup\{|x-y|/c,\,\,x\in\Gamma,\,\,y\in\Omega_0\}$.
\end{theorem}
\begin{theorem}\label{thm:uniqueness} Assume that $\partial \Omega$ is real analytic.
Then the field measurement $u|_{[0,T]\times\Gamma}$ determines the source trajectory $p$ and intensity $q$ uniquely.
\end{theorem}
\begin{theorem}\label{thm:stability}  Let $u_j $ be  the solution  of the system~\eqref{eqn:system}
with $q= q_j,$ and $p=p_j$ for $j=1, 2$.  Assume that $\|q_j\|_{C^1([0, T_0])} \leq M_q $ and
$ \|p_j\|_{C^2([0, T_0])} \leq M_p$ for  $j=1, 2$.\\

Let $\rho \in ]0, T_0[$,  $q_0$ and $c_0$ be  strictly positive constants. Assume 
\bean \label{Qassumption}
 q_j (t) \geq q_0, \quad  \forall t\in [0, T_0-\rho], 
 \eean
and  
 \bean\label{Passumption}
1- c^{-1}p_j^\prime(t)\geq c_0, \quad  \forall
t\in [0, T], \; j=1, 2. 
 \eean
 Let \begin{eqnarray*}
\varepsilon \hspace{-1mm}  = \hspace{-1mm} \sup_{\widehat  k\in S^2}  \hspace{-1mm} 
\int_{\pmb R \times \partial \Omega} \hspace{-2mm} \left(c^{-1} |\partial_t u(t-c^{-1}\widehat k \cdot x, x)|  + |\partial_\nu u(t-c^{-1}\widehat k \cdot x, x)| \right)  \hspace{-1mm} H(t-c^{-1}\widehat k \cdot x)d\sigma(x) dt, \nonumber
\end{eqnarray*}
where $u = u_1-u_2$. Then  for $\varepsilon \in ]0, q_0(T_0-\rho)[$, we have 
 \bean
 \| q_1- q_2\|_{C^0([0, T_0-\rho])} + \| p_1- p_2\|_{C^0([0, T_0-\rho])}  \leq C \varepsilon,
 \eean
 where $H$ is the Heaviside function, and   $C>0$ is a constant that  only depends on $c_0,$  $q_0,$  $T_0, $ $c, $ $M_q$, and  $M_p$.
 
 \end{theorem}

We prove Theorems~\ref{thm:regularity}--\ref{thm:stability} in Section~\ref{sec:proofs}, and describe our framework for Bayesian inference of source amplitudes and trajectories in Section~\ref{sec:Bayes}. Finally, we show and discuss our numerical results in Section~\ref{sec:numerical}. 
\section{Proofs of theorems~\ref{thm:regularity}--\ref{thm:stability}}\label{sec:proofs}
\subsection{Proof of Theorem~\ref{thm:regularity} (regularity of solution of the forward problem)}
Using the fact that
\begin{align*}
\int_{\omega\in S^2}\exp\left(-ic(t-s)|\xi|\frac{\xi}{|\xi|}\omega\right)d\omega&=\int_0^{\pi}\int_0^{2\pi}\exp\left(-i|\xi|c(t-s)\cos\theta\right)\sin\theta d\varphi d\theta\\&=\frac{4\pi}{c(t-s)}\frac{\sin(|\xi|c(t-s))}{|\xi|},
\end{align*}
we get
\begin{align}\label{eqn:uhat}
\widehat{u}(t,\xi)&=\frac{c^2}{4\pi}\int_0^t(t-s)q(s)\int_{\omega\in S^2}(\delta_{c(t-s)\omega+p(s)})_x(e^{-ix\xi})d\omega ds\nonumber\\&=\frac{c}{|\xi|}\int_0^tq(s)\sin(c|\xi|(t-s))e^{-i\xi p(s)}ds,
\end{align}
Since $u(t,\cdot)\in\mathcal{E}'^0(\pmb{R}^3)$, we know by the Paley-Wiener-Schwartz theorem~\cite[Theorem 7.3.1]{HI} that $\widehat{u}(t,\cdot)$ is entire in $\pmb{C}^3$ and that $|\widehat{u}(t,\cdot)|$ is bounded in $\pmb{R}^3$ by a constant $\widetilde{C}$; in particular, the singularity at zero in~\eqref{eqn:uhat} is removable, and for every $\xi\in\pmb{R}^3$, $t\in[0,T]$, we have
\begin{align*}
2c^{-1}|\xi||\widehat{u}(t,\xi)|&\le\left|\int_0^tq(s)\exp\left(i|\xi|f_+(s,t,\xi/|\xi|)\right)ds\right|\\&+\left|\int_0^tq(s)\exp\left(i|\xi|f_-(s,t,\xi/|\xi|)\right)ds\right|,
\end{align*}
where $f_{\pm}(s,t,\omega)=\pm c(t-s)-p(s).\omega$ for $s\in[0,T]$, $t\in[0,T]$, and $\omega\in S^2$. The real-valued functions $\partial_sf_{\pm}(s,t,\omega)=\mp c-p'(s).\omega$ are nonzero for all $s$, $t$, and $\omega$, so we can use integration by parts to write
\begin{align*}
\left|\int_0^tq(s)e^{i|\xi|f_{\pm}(s,t,\xi/|\xi|)}ds\right|&=\left|\int_0^tq(s)\left(\frac{1}{i|\xi|\partial_sf_{\pm}(s,t,\xi/|\xi|)}\frac{\partial}{\partial s}\right)^me^{i|\xi|f_{\pm}(s,t,\xi/|\xi|)}ds\right|\\&=|\xi|^{-m}\left|W(t,\xi/|\xi|)+\int_0^tQ(s,t,\xi/|\xi|)e^{i|\xi|f_{\pm}(s,t,\xi/|\xi|)}ds\right|\\&\le|\xi|^{-m}\left(|W(t,\xi/|\xi|)|+\int_0^t|Q(s,t,\xi/|\xi|)|ds\right),
\end{align*}
for some functions $W$ and $Q$ continuous w.r.t. all their arguments because $q\in C^m$ and $p\in C^{m+1}$, and because all $s$-derivatives of $f_{\pm}$ are continuous w.r.t. $t$ and $\omega$. But the sets $[0,T]$ and $S^2$ are compact, so there is a finite constant $C$, independent of $t$ and of $\xi/|\xi|$, satisfying $|\widehat{u}(t,\xi)|\le C|\xi|^{-m-1}$ for all $t\in[0,T]$, $\xi\in S^2$. Thus, for $\sigma<m-1/2$,
\begin{align*}
\|(1+|\xi|^2)^{\sigma/2}\widehat{u}(t,\cdot)\|_{L^2(\pmb{R}^3)}^2&\le 2^{\sigma}4\pi\widetilde{C}^2/3+4\pi C^2\int_{r=1}^{\infty}r^{-2m}(1+r^2)^{\sigma}dr\\&\le2^{\sigma+2}\pi \left(\widetilde{C}^2/3+C^2\left[\frac{r^{2(\sigma-m)+1}}{2(\sigma-m)+1}\right]_{r=1}^{\infty}\right)<\infty.
\end{align*}
To prove the second part of the theorem, we use the fact that~\cite[p. 130]{Duistermaat}
\begin{align*}
\singsupp{E_+(t-s)\ast(q(s)\delta_{p(s)})}&\subset\singsupp{E_+(t-s)}+\singsupp{q(s)\delta_{p(s)}}\\&=\{y+z,\,\,|y|=c(t-s),\,\,z=p(s)\}
\end{align*}
for $t\in[0,T]$ and $s\in[0,T_0]$.
Specifically,
\[
\singsupp{u(t,\cdot)}\subset\bigcup_{s=0}^{\min\{t,T_0\}}\{y+z,\,\,|y|=c(t-s),\,\,z=p(s)\},\quad t\in[0,T],
\]
and finally
\begin{align*}
\bigcup_{t=0}^T\singsupp{u(t,\cdot)}&\subset\bigcup_{t=0}^{T_0}\bigcup_{s=0}^{\min\{t,T_0\}}\{y+z,\,\,|y|=c(t-s),\,\,z=p(s)\}\\&=\bigcup_{s=0}^{T_0}\overline{B_{p(s)}(c(T_0-s))}.
\end{align*}
This completes the regularity result of the theorem. Finally, we deduce from \eqref{tau} that $\tau(t,x) $ satisfies 
\bean \label{tau2}
 t- \sup\{|x-y|/c,\,\,x\in\Gamma,\,\,y\in\Omega_0\} \leq \tau(x,t)  \leq t \quad \textrm{for all } t>0, \, x\in \Gamma.
\eean
Since $q(t)$ is supported in $[0,T_0]$, we obtain that $ u(t, \cdot)|_\Gamma$ is also supported
in  $[0,T^*]$, which finishes the proof of Theorem~\ref{thm:regularity}.

\subsection{Proof of Theorem~\ref{thm:uniqueness} (uniqueness of solution of the inverse problem)}

Under the assumptions in Section~\ref{sec:introduction}, the partial Fourier transform of $q(t)\delta_{p(t)}$ w.r.t. $t$ is well-defined in $\mathcal{E}'^0(\pmb{R}^3)$, with support in $p([0,T_0])$. Indeed, for any $\phi\in C_0^{\infty}(\pmb{R}^3)$ and any real $\omega$, we have
\begin{align*}
\left|\left(\int_{t\in\pmb{R}_+}e^{-i\omega t}q(t)\delta_{p(t)}dt\right)(\phi)\right|&=\left|\int_{t=0}^{T_0}e^{-i\omega t}q(t)\phi(p(t))dt\right|\\&\le T_0\|q\|_{L^{\infty}([0,T_0])}\|\phi\|_{L^{\infty}(p([0,T_0])\cap{\rm supp}\,\phi)}.
\end{align*}
Now assume $u_1$ and $u_2$ satisfy~\eqref{eqn:system} with sources $q_1(t)\delta_{p_1(t)}$ and $q_2(t)\delta_{p_2(t)}$, respectively, and such that $u_1|_{[0,T]\times\Gamma}\equiv u_2|_{[0,T]\times\Gamma}$. Writing $u=u_1-u_2$ and $f(t,\cdot)=q_1(t)\delta_{p_1(t)}-q_2(t)\delta_{p_2(t)}$, we have
\begin{equation}\label{eqn:system2}
\left\{
\begin{array}{rcll}
c^{-2}\partial^2_tu-\Delta u&=&f(t,\cdot)&\text{in}\,\,C_0(\pmb{R},\mathcal{E}'^0(\pmb{R}^3)),\\
u(x,0)&=&0,&x\in\pmb{R}^3,\\
\partial_tu(x,0)&=&0,&x\in\pmb{R}^3,\\
u|_{\Gamma\times\pmb{R}}&\equiv&0.
\end{array}\right.
\end{equation}
The last equation in~\eqref{eqn:system2} follows from~\eqref{eqn:solution} and the fact that $u$ is supported in $[0,T]$. Taking the partial Fourier transform of the wave equation in~\eqref{eqn:system2} with respect to $t$, we get $(\Delta+(\omega/c)^2)\widehat{u}(\omega,x)=-\widehat{f}(\omega,x)$, so, since $\supp{\widehat{f}(\omega,\cdot)}\cap\partial\Omega=\emptyset$, we have $(\Delta+(\omega/c)^2)\widehat{u}(\omega,x)=0$ for $x$ in a nonempty open neighborhood $N$ of $\partial\Omega$. The ellipticity of the Helmholtz operator $\Delta+(\omega/c)^2$ now implies that $\widehat{u}(\omega,x)$ is real-analytic w.r.t. $x$ in $N$. We have furthermore that $\widehat{u}(\omega,x)\equiv0$ for $x\in\Gamma$. If $\Gamma$ is boundaryless then $\Gamma=\partial\Omega$. Assuming $\partial\Gamma\neq\emptyset$, pick $x\in\partial\Gamma$. There is an open neighborhood $N_x\subseteq N$ of $x$ and a biholomorphism $\phi$ mapping $N_x$ to an open neighborhood of the origin in $\pmb{R}^3$, such that $\phi(x)=0$, $\phi(N_x\cap\Omega)\subset\pmb{R}^2\times\pmb{R}_-$, $\phi(N_x\cap\complement\Omega)\subset\pmb{R}^2\times\pmb{R}_+$, and $\phi(N_x\cap\partial\Omega)\subset\pmb{R}^2\times\{0\}$. The function $\widehat{U}(\omega,x')=\widehat{u}(\omega,\phi^{-1}(x',0))$ is real analytic in $\phi(N_x\cap\partial\Omega)$, and $\widehat{U}\equiv0$ in $\phi(N_x\cap\Gamma)$, which implies $\widehat{U}\equiv0$ in $\phi(N_x\cap\partial\Omega)$ by the identity theorem for holomorphic functions, as $\phi(N_x\cap\Gamma)$ is a nonempty open subset of $\phi(N_x\cap\partial\Omega)$. Since $\partial\Omega$ can be covered by a finite number of biholomorphic local charts, we conclude that $\widehat{u}(\omega,\cdot)\equiv0$ at $\partial\Omega$. Thus, for any choice of open $\Gamma\subseteq\partial\Omega$, the function $\widehat{u}(\omega,\cdot)$ satisfies the exterior Helmholtz problem
\begin{equation}\label{eqn:sysH}
\left\{\begin{array}{rcll}
(\Delta+(\omega/c)^2)\widehat{u}(\omega,x)&=&0,\quad&x\in\pmb{R}^3\setminus\overline{\Omega},\\
\widehat{u}(\omega,x)&=&0,\quad&x\in\partial\Omega,\\
\lim_{|x|\rightarrow\infty}|x|(\partial_{|x|}\widehat{u}-i(\omega/c)\widehat{u})&=&0,\quad&\text{uniformly in}\,\,x/|x|\in S^2.
\end{array}\right. 
\end{equation} 
 The Sommerfeld radiation condition in $\pmb{R}^3$, which is the last equation in~\eqref{eqn:sysH},  is satisfied by $\widehat{u}$ because the support of the source $\widehat{f}(\omega,\cdot)$ has a compact support; indeed, $\widehat{f}(\omega,x)=0$ if $x\notin p_1([0,T_0])\cup p_2([0,T_0])$. But~\eqref{eqn:sysH} is known~\cite[Theorem 2.6.5, p. 102]{Nedelec} to have only the trivial solution.
In particular $\widehat{u}= \partial_{\nu}\widehat{u}=0$ at $\partial\Omega$, and hence
$\widehat{u}$  solves the interior Helmholtz problem 
 \begin{equation}\label{eqn:sysHinte}
\left\{\begin{array}{rcll}
 (\Delta+(\omega/c)^2)\widehat{u}(\omega, x)&=&-\widehat{f}(\omega, x),\quad&\text{in}\,\,\Omega,\\
\widehat{u}(\omega, w)&=&0,\quad&x\in\partial\Omega,\\
\partial_{\nu}\widehat{u}(\omega, x)&=&0,\quad&x\in\partial\Omega.
\end{array}\right.
\end{equation}
Since  $\widehat{f}(\omega,\cdot)$ has a compact support $ p_1([0,T_0])\cup p_2([0,T_0])$,  which has
 a zero  two-dimensional Hausdorff measure,  we deduce from  the uniqueness of Cauchy problem,  $\widehat{u} = 0 $  a. e. in $\Omega$. On the other hand  a simple calculation shows that $\widehat{f}(\omega,\cdot) 
\in W^{-1, p}(\pmb{R}^3)$ for all $p \in [1, \frac{3}{2}[$, which regarding to the elliptic regularity 
 implies  that  $\widehat{u}(\omega, \cdot)  \in W^{1, p}_{loc}(\pmb{R}^3) $ for all $ p\in [1, \frac{3}{2}[$.
 Consequently   $\widehat{u} = 0 $   in $\mathcal{E}'^0(\pmb{R}^3)$. We  then obtain $\widehat{f}(\omega,\cdot) = 0$,  which  completes the proof of Theorem~\ref{thm:uniqueness}.
\subsection{Proof of Theorem~\ref{thm:stability} (stability of solution of the inverse problem)}

Let $u_1$ and $u_2$ satisfy~\eqref{eqn:system} with sources  $q_1(t)\delta_{p_1(t)}$ and $q_2(t)\delta_{p_2(t)}$, respectively, and such that $u_1-u_2 $ and  $\partial_\nu (u_1 -u_2)$
 are  given on $[0,T]\times\partial \Omega$. Writing $u=u_1-u_2$ and $f(t,\cdot)=q_1(t)\delta_{p_1(t)}-q_2(t)\delta_{p_2(t)}$, we have
\begin{equation}\label{eqn:system22}
\left\{
\begin{array}{rcll}
c^{-2}\partial^2_tu-\Delta u&=&f(t,\cdot)&\text{in}\,\,C_0(\pmb{R},\mathcal{E}'^0(\pmb{R}^3)),\\
u(x,0)&=&0,&x\in\pmb{R}^3,\\
\partial_tu(x,0)&=&0,&x\in\pmb{R}^3,
\end{array}\right.
\end{equation}
Applying  the partial Fourier transform to  the wave equation in~\eqref{eqn:system22} with respect to $t$, we obtain 
\begin{equation}\label{eqn:sysH2}
\left\{\begin{array}{rcll}
(\Delta+(\omega/c)^2)\widehat{u}(\omega,x)&=&-\widehat{f}(\omega,x),\quad&x\in\pmb{R}^3,\\
\lim_{|x|\rightarrow\infty}|x|(\partial_{|x|}\widehat{u}-i(\omega/c)\widehat{u})&=&0\quad&\text{uniformly in}\,\,x/|x|\in S^2.
\end{array}\right.
\end{equation}
Let $\widehat{k}\in S^2$ be fixed. 
Multiplying the Helmholtz equation by $e^{-i\frac{\omega}{c}\widehat k \cdot x}$, and integrating by parts
yield
\begin{eqnarray*}
\int_{t\in\pmb{R}_+}q_1(t)e^{-i(\omega/c)\widehat{k}.p_1(t)}e^{-i\omega t} dt-\int_{t\in\pmb{R}_+} q_2(t)e^{-i(\omega/c)\widehat{k}.p_2(t)}e^{-i\omega t}dt\\
 = -\int_{\partial \Omega} \left(i \frac{\omega}{c}\widehat k \cdot \nu \widehat u  + \partial_\nu \widehat u\right) e^{-i\frac{\omega}{c}\widehat k \cdot x} d\sigma(x). 
\nonumber
\end{eqnarray*}
Making the change of variables $t = \varphi_j(\tau)$  where $\tau=\varphi_j(\tau)+c^{-1}\widehat{k}.p_j(\varphi_j(\tau))$, $j=1,2$, respectively in the two 
integrals  on the left side give
\begin{eqnarray}\label{eqn:00}
\int_{\tau\in\pmb{R}}\left(\frac{\chi_1(\tau)q_1(\varphi_1(\tau))}{c^{-1}\widehat{k}.p_1'(\varphi_1(\tau))+1}-\frac{\chi_2(\tau)q_2(\varphi_2(\tau))}{c^{-1}\widehat{k}.p_2'(\varphi_2(\tau))+1}\right)e^{-i\omega\tau}d\tau = \\ 
-\int_{\partial \Omega} \left(i \frac{\omega}{c}\widehat k \cdot \nu \widehat u + \partial_\nu \widehat u\right) e^{-i\frac{\omega}{c}\widehat k \cdot x} d\sigma(x),
\nonumber
\end{eqnarray}
where $\chi_{j}(\tau) = H( \tau- \tau_j) $, with
$\tau_j = \varphi_j^{-1}(0)= c^{-1}\widehat{k}.p_j(0), j=1,2 $, and $H$ is the Heaviside function. Since $\varphi_j,  j=1,2, $
are strictly increasing functions, we also have $\chi_{j}(\tau) = H(\varphi_j(\tau)),  j=1,2 $.   \\

Applying the Fourier transform inverse both sides, we get 
\begin{eqnarray}\label{eqn:77}
\frac{\chi_1(\tau)q_1(\varphi_1(\tau))}{c^{-1}\widehat{k}.p_1'(\varphi_1(\tau))+1}-\frac{\chi_2(\tau)q_2(\varphi_2(\tau))}{c^{-1}\widehat{k}.p_2'(\varphi_2(\tau))+1} = \\ 
-\int_{\partial \Omega} \left(c^{-1}\widehat k \cdot \nu \partial_t u(\tau-c^{-1}\widehat k \cdot x, x) + \partial_\nu u(\tau-c^{-1}\widehat k \cdot x, x) \right) H(\tau-c^{-1}\widehat k \cdot x)d\sigma(x), \nonumber
\end{eqnarray}
for all $\tau \in \pmb R$. \\

Recall that Theorem \ref{thm:regularity}  implies that   $u\in C^{\infty}(]0, T[, H^{m-1-\varepsilon}(\partial \Omega))$, and hence the right hand side term is indeed well defined. \\

Notice that $\varphi_j(\tau), \tau_j,  j=1, 2, $  as well as $\tau_j, j=1, 2 $, are also functions 
of $\widehat k$.  Next our strategy  is to first estimate $\| p_1(0) - p_2(0)\| $ in terms of the boundary
measurements. 

Without loss of generality one can assume that $\tau_1 \leq \tau_2$. Integrating \eqref{eqn:77} both sides over
$]\tau_1, \tau_2[$, yields 
\begin{eqnarray}\label{eqn:505}
\int_{0}^{\varphi_1(\tau_2)}q_1(s) ds = \int_{\tau_1}^{\tau_2}\frac{q_1(\varphi_1(t))}{c^{-1}\widehat{k}.p_1'(\varphi_1(t))+1} dt \leq \varepsilon. 
\end{eqnarray}
We claim that $ \varphi_1(\tau_2) \in ]0, T_0-\rho[$. Indeed, assuming that $\varphi_1(\tau_2)>
  T_0-\rho$, we deduce from  \eqref{Qassumption}, and \eqref{eqn:505},  $ q_0(T_0-\rho) < \varepsilon$, which is in contradiction with the fact that   $\varepsilon \in ]0, q_0(T_0-\rho)[$. \\
  
Therefore $ \varphi_1(\tau_2) \in ]0, T_0-\rho[$, and hence 
   
\bean \label{firstestimate}
|\tau_2 -\tau_1| \leq 2 q_0^{-1}\varepsilon. 
\eean
Since the right hand side term  is independent of $\widehat k$, we also get
 
\bean \label{firstmainestimate}
\| p_1(0) - p_2(0)\|  \leq 2 c q_0^{-1}\varepsilon. 
\eean


 
 Let $\widehat k_0= \widehat k(t_0)\in S^2$ such that $\widehat k_0\cdot(p_1(t_0)-p_2(t_0)) = 0$.  
We deduce from \eqref{eqn:77} that 
\begin{equation}\label{eqn:equa2}
\left| \frac{q_1(t_0)}{\pm c^{-1}\widehat{k}_0.p_1'(t_0)+1}-\frac{q_2(t_0)}{\pm c^{-1}\widehat{k}_0.p_2'(t_0)+1}\right| \leq \varepsilon.
\end{equation}

Simple calculations yield 
\bea
|q_1(t_0)- q_2(t_0)| \leq 4\varepsilon.
\eea

Since $t_0$ is arbitrarily, we obtain 
\bean  \label{secondmainestimate}
\|q_1- q_2\|_{C^0([0, T_0])} \leq 4\varepsilon. 
\eean

Integrating the identity \eqref{eqn:77} both sides over $]-\infty, r[$, with $r>0$, gives 

\begin{eqnarray}\label{eqn:808}
\left| \int_{0}^{\varphi_1(r)}q_1(s) ds - \int_{0}^{\varphi_2(r)}q_2(s) ds \right|  \leq \varepsilon, \quad 
\end{eqnarray}
Combining  \eqref{secondmainestimate}  and  \eqref{eqn:808}  lead to

\begin{eqnarray}\label{eqn:809}
\left| \int^{\varphi_2(r)}_{\varphi_1(r)}q_2(s) ds \right|  \leq(1+4T_0) \varepsilon,
\end{eqnarray}
for all $r\in ]r_1, r_2[$ with 
\[
 r_2= \min_{j=1,2} \varphi_j^{-1}(T_0-\rho), \; \;
r_1= \max_{j=1,2}\varphi_j^{-1}(0)= \max_{j=1,2}\tau_j.
\]
Notice that $r\in ]r_1, r_2[$ is equivalent to $\varphi_j(r) \in ]0, T_0-\rho[, j=1, 2$. \\

 Using assumption \eqref{Qassumption},
we find 
\begin{eqnarray}\label{thirdmainestimate}
\left|\varphi_2(r)- \varphi_1(r) \right|  \leq q_0^{-1}(1+4T_0) \varepsilon,
\end{eqnarray}
for all $r\in  ]r_1, r_2[.$

Without loss of generality we  further assume  $r_2= \varphi_1^{-1}(T_0-\rho).$\\

Further $C$ designate a strictly positive constant  that  only depends on $c_0,$  $q_0,$  $T_0, $ $c, $ $M_q,$ and $M_p.$\\

Combining identities \eqref{eqn:77}, \eqref{secondmainestimate}, and \eqref{thirdmainestimate}, we  obtain 

\bea 
 |q_1(\varphi_1(r))|\left|\widehat k \cdot (p_1'(\varphi_1(r))- p_2'(\varphi_1(r))) \right| \leq C\varepsilon, \quad \forall
 r\in  ]r_1, r_2[, \, \forall \widehat k\in S^{2},
\eea
which in turn implies 

\bean \label{eqn:909}
 \left\|p_1'(t)- p_2'(t) \right\|_{C^0([\varphi_1(r_1),  T_0-\rho])} \leq   C\varepsilon.
 \eean
We deduce from \eqref{firstestimate} the following inequality 
\bean \label{eq_inter}
|r_1-\tau_1 | \leq 2 q_0^{-1}\varepsilon. 
\eean

Using estimates \eqref{firstmainestimate}, \eqref{eqn:909}, and \eqref{eq_inter}, we finally get

\bean \label{fourthmainestimate}
 \left\|p_1(t)- p_2(t) \right\|_{C^0([0, T_0-\rho])} \leq   C\varepsilon,
 \eean
which achieves the proof of the theorem.



\section{Bayesian inference of source trajectory and intensity}\label{sec:Bayes}

We next describe our setup for the Bayesian inference of the source trajectory and intensity. We confine the source trajectory $p(t)$ to the $xy$-plane, setting $p(t)=(p_x(t),p_y(t),0)$, and impose GP (Gaussian process) priors on $p_x(t)$, $p_y(t)$, and $q(t)$. In our case, a GP is a stochastic process $g:\pmb{R}\rightarrow\pmb{R}$ such that for any $d\in\pmb{N}$ and any $t_1,\dots,t_d\in[0,T_0]$ the distribution of the $d$-tuple $(g(t_1),\dots,g(t_d))$ is multivariate Gaussian. We write $g\sim\mathcal{GP}(m,k)$, where $m(t)$, $t\in[0,T_0]$, is the mean function and $k(t,t')$, $t,t'\in[0,T_0]$, is the covariance function. 
The mean and the covariance functions influence the dimensionality, smoothness, stationarity, periodicity, and other properties of the realizations of the GP. We choose the squared-exponential (SE) kernel
\[
k_{\rm SE}(t,t')=\kappa^2\exp(-(t-t')^2/2\ell^2),\quad t,t'\in[0,T],
\]
where the hyperparameters $\kappa$ and $\ell$ are the magnitude and the correlation length, respectively. The magnitude controls the extent to which realizations of the GP can deviate from the mean, while the correlation length determines the speed with which these realizations can oscillate (larger $\ell$ gives slower oscillation). The SE kernel results in a smooth prior on the functions sampled from the GP, and corresponds to the use of radial basis functions $\phi_j$. Also, $g\sim\mathcal{GP}(0,k_{\rm SE})$ is second-order stationary, since $k_{\rm SE}$ is isotropic. We impose independent GP priors on the latent functions $p_x(t)$, $p_y(t)$ and $q(t)$: $p_x(t)\sim\mathcal{GP}(0,k_{\rm SE}^p)$, $p_y(t)\sim\mathcal{GP}(0,k_{\rm SE}^p)$, $q(t)\sim\mathcal{GP}(0,k_{\rm SE}^q)$. The superscripts $p$ and $q$ indicate possibly separate choices of the hyperparameter values for the GP priors. Our measurement data $(X_i, U_i)_{i=1}^N$ consist of the field sampling times and sensor locations $X_i=(t_i,x_i,y_i,z_i)$ and sampled field values $U_i=u(t_i,x_i, y_i, z_i)$, where $N=N_sN_t$ is the product of the number $N_s$ of sensors and the number $N_t$ of samplings. We want to estimate $f = (p_x(t), p_y(t), q(t))$, having the forward mapping $G$ from parameter space to data space given by
\[
G(f)(X_i)=\begin{cases}\frac{c}{4\pi}\frac{q(\tau(X_i))}{|(x_i,y_i,z_i)-p(\tau(X_i))|h((x_i,y_i,z_i),\tau(X_i))},&|(x_i,y_i,z_i)-p(0)|<ct_i,\\0&{\rm otherwise.}\end{cases}
\label{eq:forward_op}
\]
The output is related to the input by $U_i = G(f)(X_i) + e$, with $e \sim \mathcal{N}(0,\beta^{-1}I)$ for some likelihood precision $\beta$. Thus, the field measurements are distributed according to
\[
\pmb{U}|{p_x, p_y, q} \sim \mathcal{N}(G(p_x, p_y, q)(\pmb{X}), \beta^{-1}\pmb{I}_{N\times N}),
\]
and the Bayesian posterior for the latent functions is given by
\[
\mathcal{P}(p_x,p_y,q|\pmb{U})=\frac{\mathcal{P}(\pmb{U}|p_x,p_y,q)\mathcal{P}(p_x)\mathcal{P}(p_y)\mathcal{P}(q)}{\mathcal{P}(\pmb{U})}.
\]
The nonlinearity of $G$ makes an analytic treatment of the posterior intractable, and we resort to a Markov Chain Monte Carlo (MCMC) numerical procedure. We evaluate the unknown functions at a set of sampling time points $\pmb \tau\in[0,T_0]^{M}$. We let $\vec{p}_x, \vec{p}_y$ and $\vec{q}$ be the unknown functions and $\vec{K}$ be the kernel evaluated at the grid $\boldsymbol{\tau}$. Note that the number of sampling points $M$, used to numerically sample the latent functions, is not related to the number of measurement point $N$. The posterior density for the latent functions sampled on $\pmb \tau$ is then
\renewcommand{\P}{\mathcal{P}}
\renewcommand{\vec}[1]{\mathbf{#1}}
\begin{align}
	\P(\vec{f}| \vec{U}) = \frac{\P(\vec{U} | \vec{f}) \P(\vec{f})}{\P(\vec{U})} \propto \P_{\N}(\vec{U}|G(\vec{f})(\vec{X}), \beta^{-1}\vec{I}) \P_{\N}(\vec{p}|0, \vec{K}_p) \P_{\N}(\vec{p}_y|0, \vec{K}_p) \P_{\N}(\vec{q}|0, \vec{K}_q)
\end{align}
and the log-likehood is
\[
\mathcal{P}(\pmb{U}|\pmb{f})=-\frac{1}{2}\beta\|\pmb{U}-G(\pmb{f})\|^2.
\]
To sample from our GP efficiently over the grid $\pmb{t}$, we use Cholesky factorization. Indeed, if $w(t)\sim\mathcal{GP}(m,k)$ then $w(\pmb{t})\sim\mathcal{N}(\pmb{m},\pmb{k})$, and if $\pmb{k}=LL^T$ is the Cholesky factorization of $\pmb{k}$ then we can sample $f(\pmb{t})=\pmb{m}+L\pmb{s}\sim\mathcal{N}(\pmb{m},\pmb{k})$ by simply sampling $\pmb{s}\sim\mathcal{N}(0,I)$. The factorization of $\pmb{k}$ need only be performed once. We use the pre-conditioned Crank-Nicholson scheme (pCN-MCMC) to select the next sample in the chain, as described in Algorithm~\ref{alg:pCN-MCMC}. The pCN-MCMC is similar to, but more efficient than Metropolis-Hastings. The main differences of pCN-MCMC relative to Metropolis-Hastings are that the proposal distribution is identical to the prior distribution $\Pi$, the proposal function is a mixture of the previous step in the chain and the new sample, and the acceptance probability is computed using only the log-likelihood instead of the log-joint density.
\begin{algorithm}
\caption{The pre-conditioned Crank-Nicholson scheme.}\label{alg:pCN-MCMC}
\begin{algorithmic}
\State assume $f\sim\Pi=\mathcal{GP}(m,k)$
\State make a starting guess $f^{(0)}\sim\Pi$
\For{$k=1$ to $K$} 
    \State draw a sample from $\psi\sim\Pi$
    \State $f^{\ast}\gets \sqrt{1-2\delta}f^{(k-1)}+\sqrt{2\delta}\psi$ \Comment{define new proposal}
    \State $A_k\gets\min\{1,\exp(L(f^{\ast})-L(f^{(k-1)}))\}$ \Comment{compute acceptance probability} \State \Comment{$L(f)$ is the log-likelihood}
    \State draw $u_k\sim \mathcal{U}(0,1)$
    \If{$u_k<A_k$}
        \State $f^{(k)}\gets f^{\ast}$ \Comment{accept new proposal}
    \Else
        \State $f^{(k)}\gets f^{(k-1)}$ \Comment{reject new proposal}
    \EndIf
\EndFor
\end{algorithmic}
\end{algorithm}
In situations where computation of the forward map is more expensive than the case we consider here, a number of methods have been suggested such as using local approximations of the forward map \cite{conrad2016accelerating}, using neural networks \cite{antil2021novel} or exploiting geometric properties of the posterior \cite{beskos2017geometric}.

\subsection{Assessing convergence}

We here discuss methods to assess the convergence of the MCMC algorithm to a stationary distribution. Inference from samples generally suffer from two main issues~\cite{gelman2013bayesian}. The first is insufficient simulation length resulting in samples that do not accurately reflect the underlying distribution. The second is correlations between the samples which reduces the effective number of samples. The first issue can be monitored by having multiple initial guesses and confirm that they converge to the same distribution. The second issue can be monitored by calculating an effective sample size for the chain to obtain the equivalent number of i.i.d samples. In practice, an effective sample size of 100 is often enough to obtain accurate posterior estimates~\cite{gelman2013bayesian}.

\subsection{Exploiting prior knowledge through conditioning} \label{sec:prior_knowledge}

One key strength of the Bayesian approach we consider here, is the ability to condition the priors given prior knowledge of the system behavior. The most common case is when the latent functions have known values at a number of points. If a function evaluated at the grid $\tau$ has the distribution $f(\pmb{\tau})\sim \mathcal{N}(\pmb{m}, \pmb{K})$ then the conditional distribution on the set of points $(\pmb{\tau}_c, f(\pmb{\tau}_c))$ is

\[
f(\pmb{\tau})|\pmb{\tau}_c,f(\pmb{\tau}_c)\sim\mathcal{N}(\pmb{m}+\pmb{k}(\pmb{K}_c+\sigma^2\pmb{I})^{-1}(f(\pmb{\tau}_c)-\pmb{m}_c),\pmb{K}-\pmb{k}(\pmb{K}_c+\sigma^2\pmb{I})^{-1}\pmb{k}^T),
\]
where $\pmb{K}_c$ is the covariance function evaluated on $\pmb{\tau}_c, \pmb{\tau}_c$, $\pmb{K}_c$ is the covariance function evaluated on $\pmb{\tau}, \pmb{\tau}_c$ and $\pmb{m}_c$ is the mean function evaluted at $\pmb{\tau}_c$.

In many cases we cannot condition on specific function values, but instead have a constraint expressed through a linear functional, i.e. $L[f] = l$. An example that we consider later is when the trajectory is known to be closed, i.e. $L[p] = p(T_0) - p(0) = 0$. If $f \sim GP(m, k)$ then the distribution of $f$ conditioned on $L[f] = l$ is $f \sim GP(m_{f|l}, k_{f|l})$ where
\begin{align}
    m_{f|l}(t) &= m(t) + \frac{L[k(\cdot, t)]}{L^2[k]}(l- L[m])\\
    k_{f|l}(t,t') &= k(t,t') - \frac{L[k(\cdot, t')]L[k(t,\cdot)]}{L^2[k]},
\end{align}
where $L^2[k]$ denotes the application of $L$ to both arguments of $k$. The ability to flexibly incorporate prior knowledge in a non-parametric way through conditioning is a key advantage of this method.

\subsection{Evaluating the forward map efficiently}

Since the numerical procedure relies on iteratively computing the forward operator in Eq. \eqref{eq:forward_op}, we need to perform this computation efficiently. Instead of solving the equation $t  = \tau +  \norm{x - p(\tau)}/c$ for the emission time $\tau$ for every measurement point $X_i$, we consider the set of emission times $\{\tau_i\}_{i =1 }^{N_\tau}$ on which we sample the latent functions. We then calculate the corresponding observation times  $t  = \tau +  \norm{x - p(\tau)}/c$ for each sensor position. These observation times differ from the measurement times $\{t_i\}_{i = 1}^{N_t}$ and we use simple linear interpolation to obtain the forward solution at the given measurement times.

\section{Numerical results}\label{sec:numerical}
\begin{figure}[ht]
\centering
\begin{minipage}[c]{0.48\textwidth}
\centering
    \includegraphics[width=\linewidth]{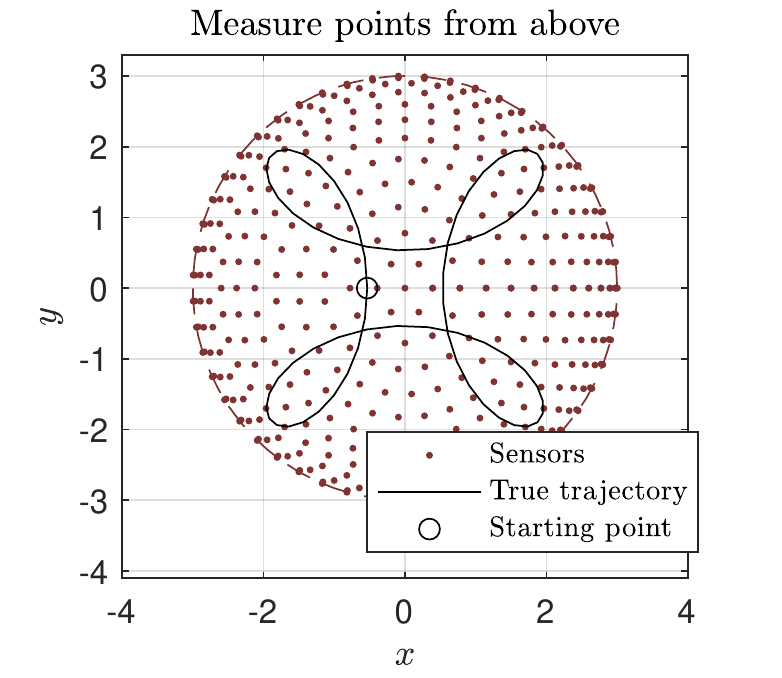}
    \end{minipage}
    \begin{minipage}[c]{0.48\textwidth}
\centering
        \includegraphics[width=\linewidth]{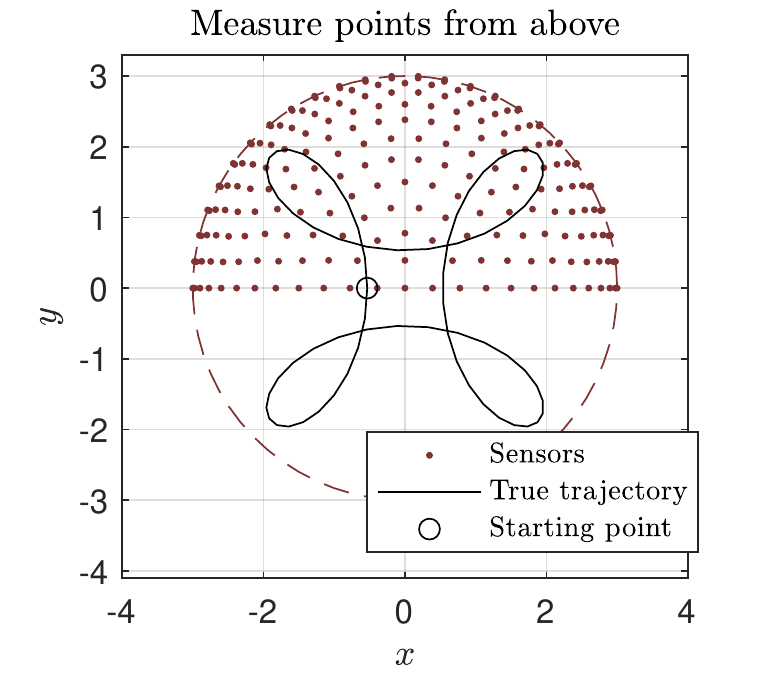}
\end{minipage}
\caption{Sensor distribution in the two setups (top view).}
    \label{fig:meas_setup} 
\end{figure}

\begin{figure}
\centering
\begin{minipage}[c]{0.49\textwidth}
\centering
    \includegraphics[width=\linewidth]{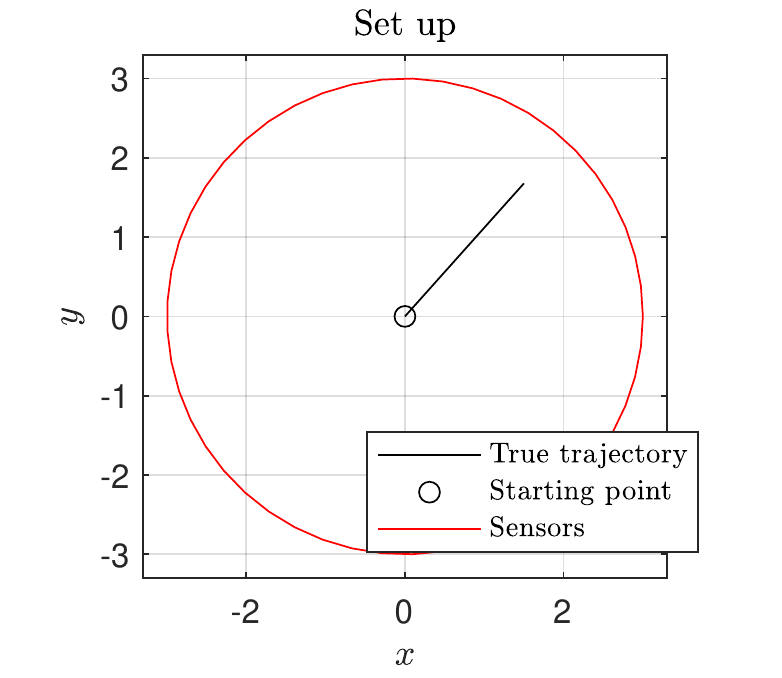}
        \includegraphics[width=\linewidth]{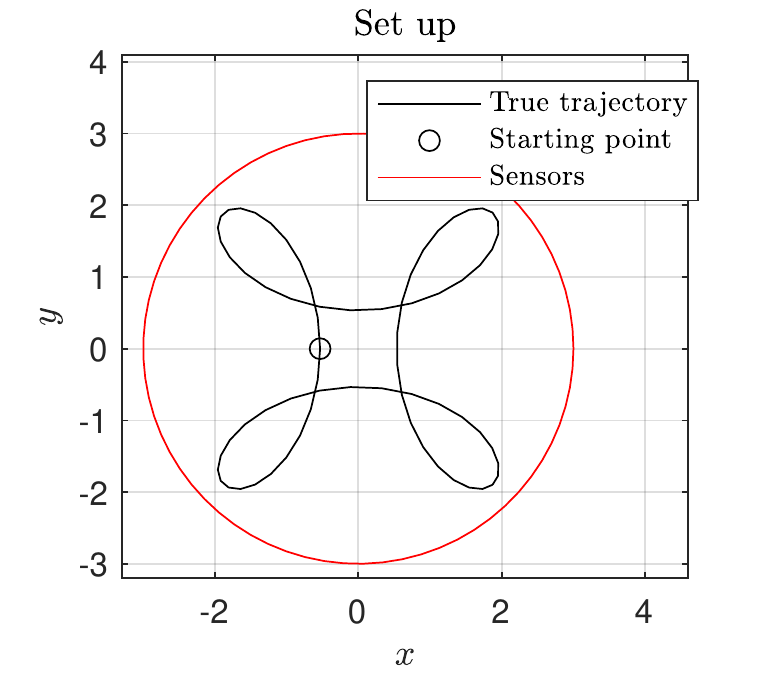}
\end{minipage}
\begin{minipage}[c]{0.49\textwidth}
\centering
    \includegraphics[width=\linewidth]{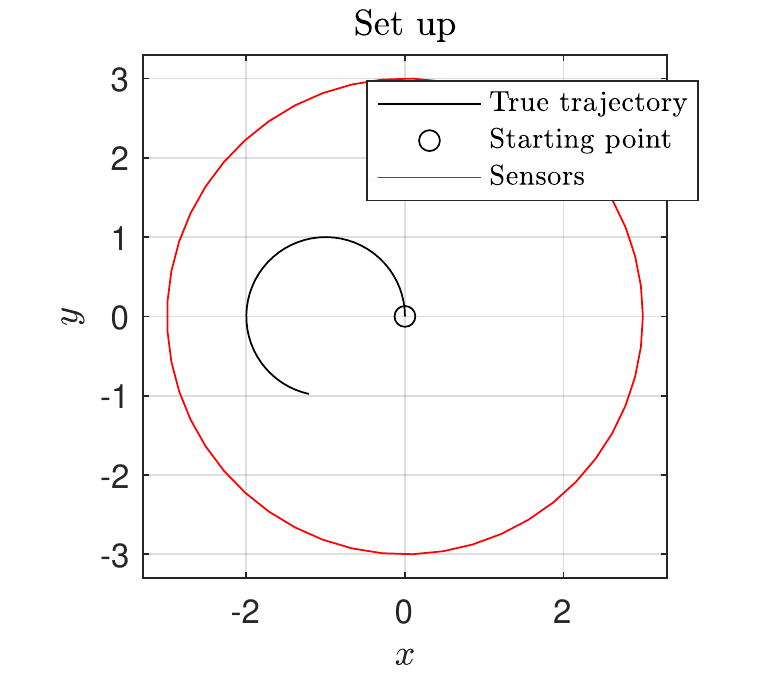}
        \includegraphics[width=\linewidth]{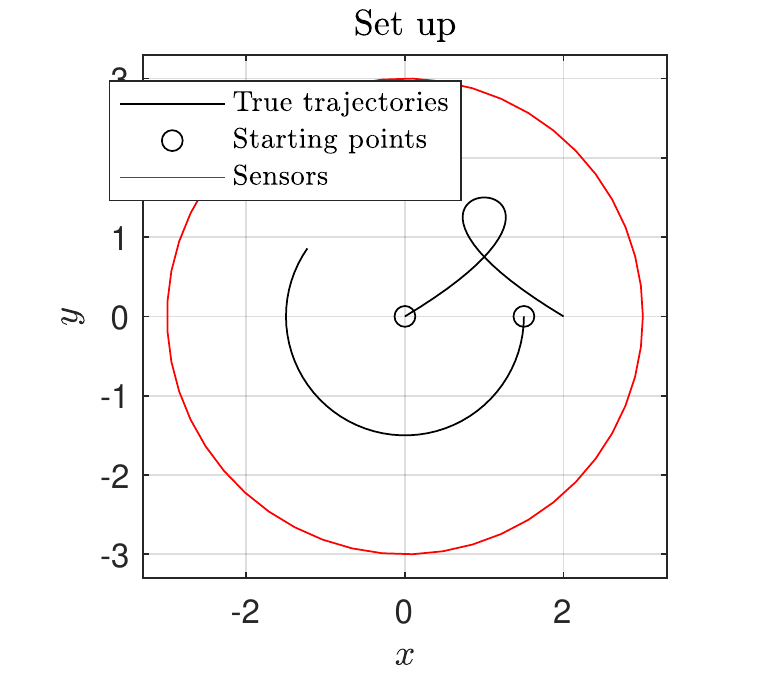}
\end{minipage}
\caption{Four different trajectories for reconstruction.}
    \label{fig:traj} 
\end{figure}

\begin{figure}
\centering
\begin{minipage}[c]{0.48\textwidth}
\centering
        \includegraphics[width=\linewidth]{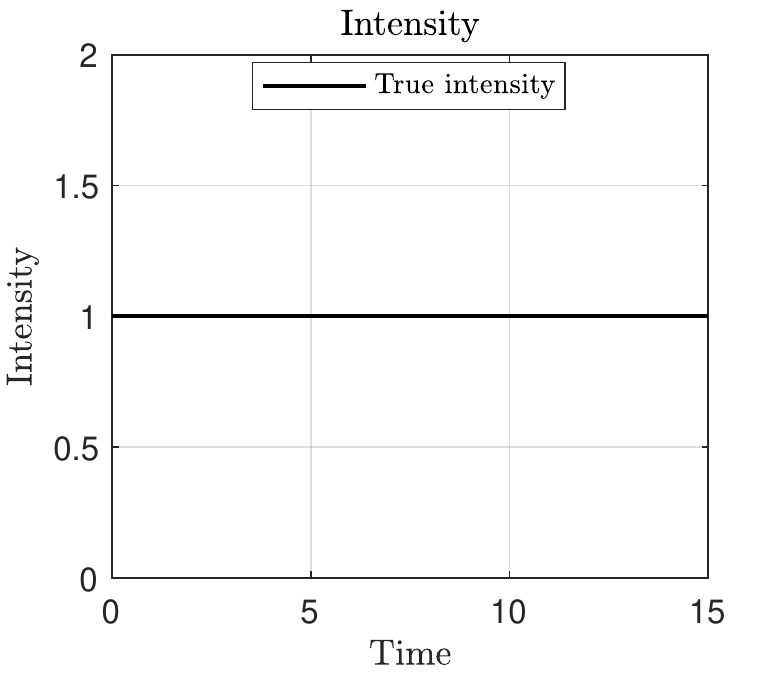}
\end{minipage}
\begin{minipage}[c]{0.48\textwidth}
\centering
        \includegraphics[width=\linewidth]{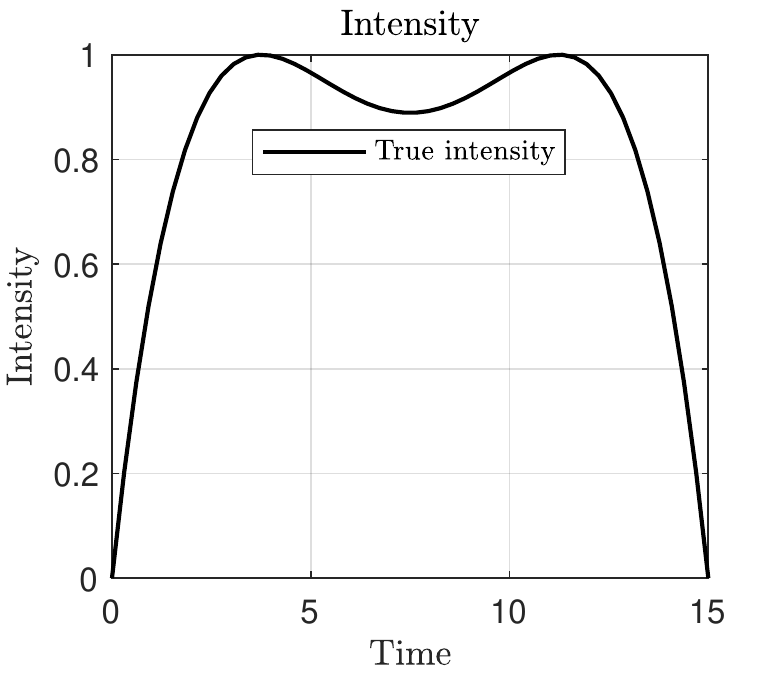}
\end{minipage}
\caption{Two different intensities for reconstruction.}
    \label{fig:inten} 
\end{figure}

We used two different measurement setups, illustrated in Figure~\ref{fig:meas_setup}, to produce the numerical results shown in this section: one with 424 sensors (field sampling points) distributed approximately uniformly over the hemisphere $\{(x,y,z)\in\mathbf{R}^3,\,\,z>0,\,\,x^2+y^2+z^2=3^2\}$, and one with 213 sensors distributed approximately uniformly over the quarter-sphere $\{(x,y,z)\in\mathbf{R}^3,\,\,z>0,\,\,y>0,\,\,x^2+y^2+z^2=3^2\}$. For visualisation purposes we here confine all source trajectories to the $xy$-plane. Consequently, identical data would be measured at the upper and lower hemispheres.
For the target sources, we chose the four trajectories shown in Figure~\ref{fig:traj}, one of which includes two distinct point sources. As intensity profiles we use the two intensities shown in Figure~\ref{fig:inten}.
All reconstructions used 6 chains of 100.000 samples from the posterior distribution generated by pCN-MCMC, each with a different initial guess of the trajectory and the intensity. The likelihood precision was set to $\beta=100$, and the pCN parameter to $\delta\in[0.001,0.0025]$, giving an acceptance ratio close to 25\%. The first 50.000 samples in each chain were discarded when calculating the posterior mean. The measuring time was $t\in[0,T]=[0, 20]$, with the source emitting during the interval $t\in[0,T_0]=[0, 15]$. For the 'long trajectory' (the closed curve in Figure~\ref{fig:traj}) $t\in[0,T]=[0, 40]$ and the source emitted during $t\in[0,T_0]=[0, 35]$.

We quantify the errors in the numerical reconstructions $\widehat{p}(t)$ and $\widehat{q}(t)$ of the posterior trajectory and intensity, respectively, as follows:
\begin{subequations}
\label{eq:errors}
\begin{align}
\text{trajectory error}&=\sqrt{\frac{1}{T_0}\int_0^{T_0}\|\widehat{p}(t)-p(t)\|^2dt} \\&\approx\sqrt{\frac{1}{N_t}\sum_{n=1}^{N_t}\|\widehat{p}_n-p_n\|^2},\notag
\end{align}
\begin{align}
\text{intensity error}&=\sqrt{\frac{1}{T_0}\int_0^{T_0}(\widehat{q}(t)-q(t))^2dt}\\&\approx\sqrt{\frac{1}{N_t}\sum_{n=1}^{N_t}(\widehat{q}_n-q_n)^2}.\notag
\end{align}
\end{subequations}

\subsection{Case 1} 
\begin{figure}[h]
\centering
\includegraphics[width=0.5\linewidth]{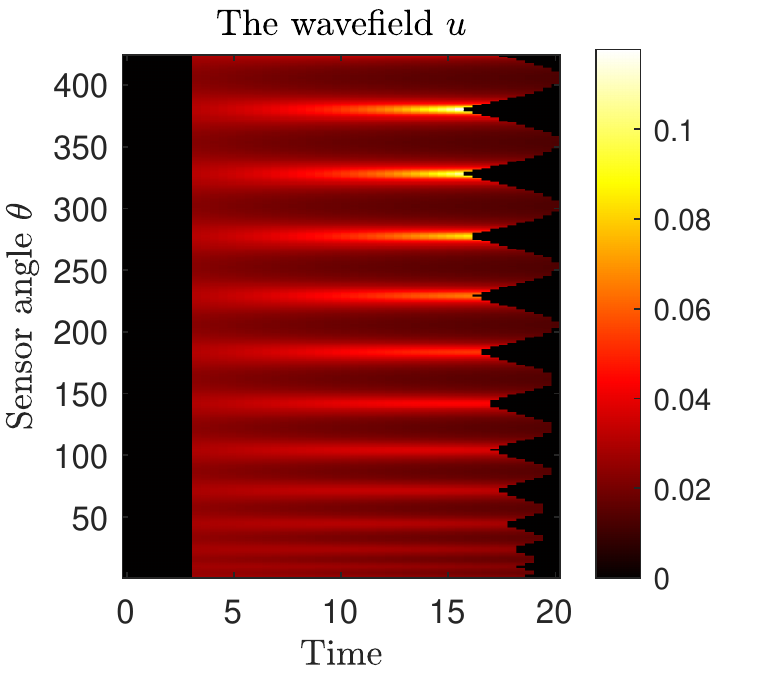}
\caption{The measured wave field in Case 1.}\label{fig:case1_measurements}
\end{figure}

\begin{figure}
\centering
\begin{minipage}[c]{0.32\textwidth}
\centering
    \includegraphics[width=\linewidth]{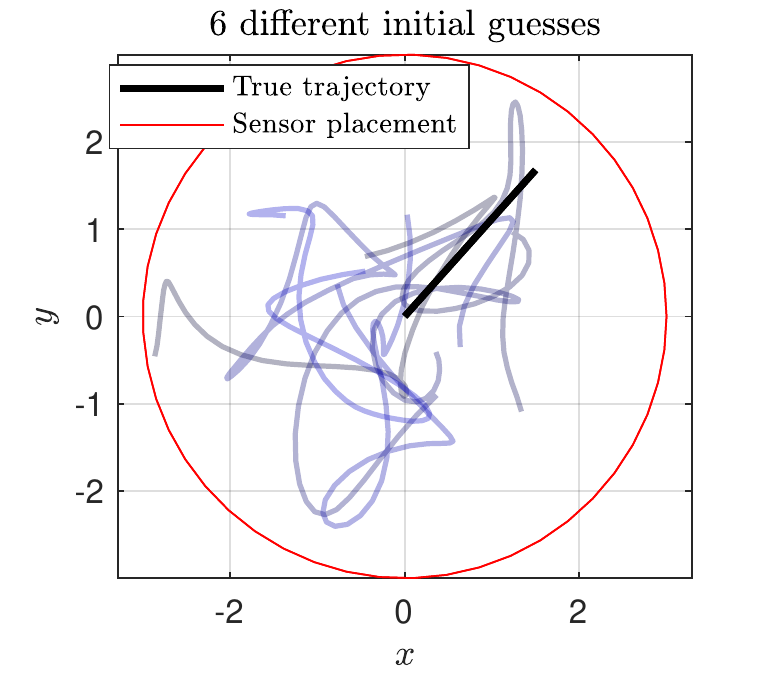}
    $\ell = 2$, $\kappa = 1$
\end{minipage}
\begin{minipage}[c]{0.32\textwidth}
\centering
    \includegraphics[width=\linewidth]{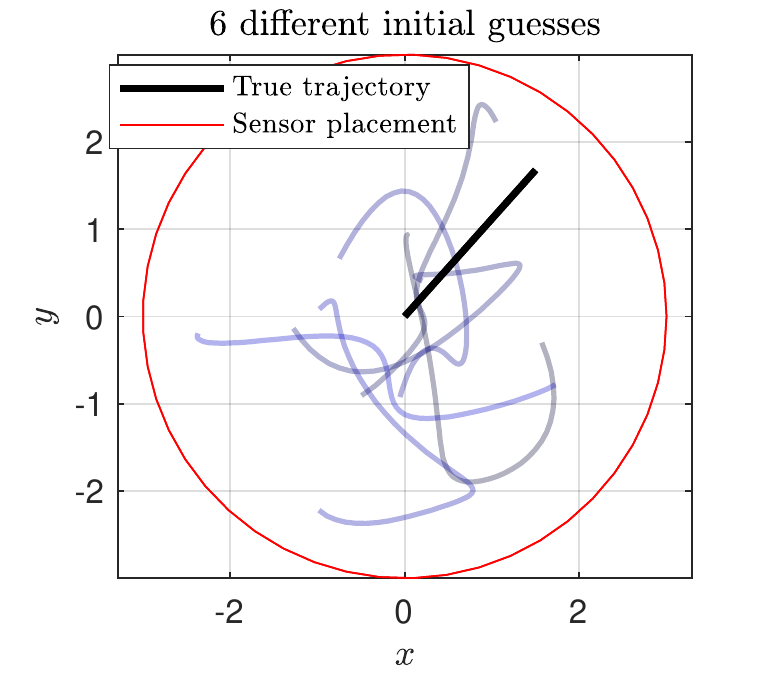}
     $\ell = 4$, $\kappa = 1$
\end{minipage}
\begin{minipage}[c]{0.32\textwidth}
\centering
    \includegraphics[width=\linewidth]{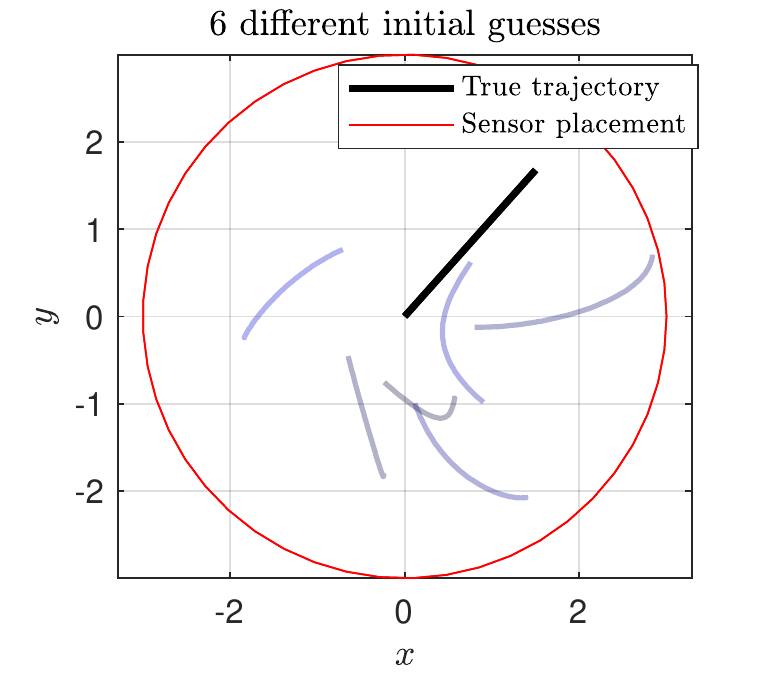}
    $\ell = 15$, $\kappa = 1$
\end{minipage}
\caption{Three types of initial guesses.}
    \label{fig:case1_initial} 
\end{figure}

\begin{figure}[h]
\begin{minipage}[c]{0.32\textwidth}
\centering
    \includegraphics[width=\linewidth]{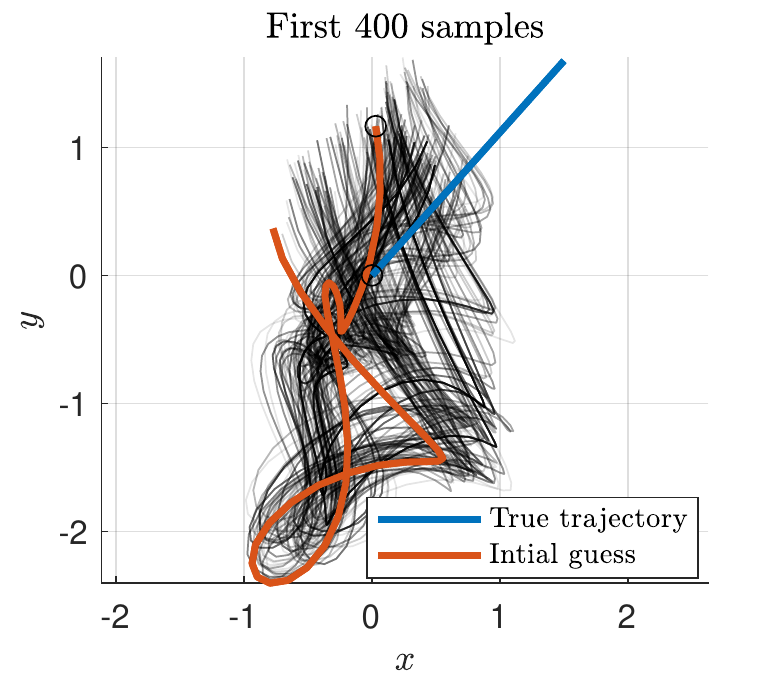}
$\ell = 2$, $\kappa = 1$
\end{minipage}
\begin{minipage}[c]{0.32\textwidth}
\centering
    \includegraphics[width=\linewidth]{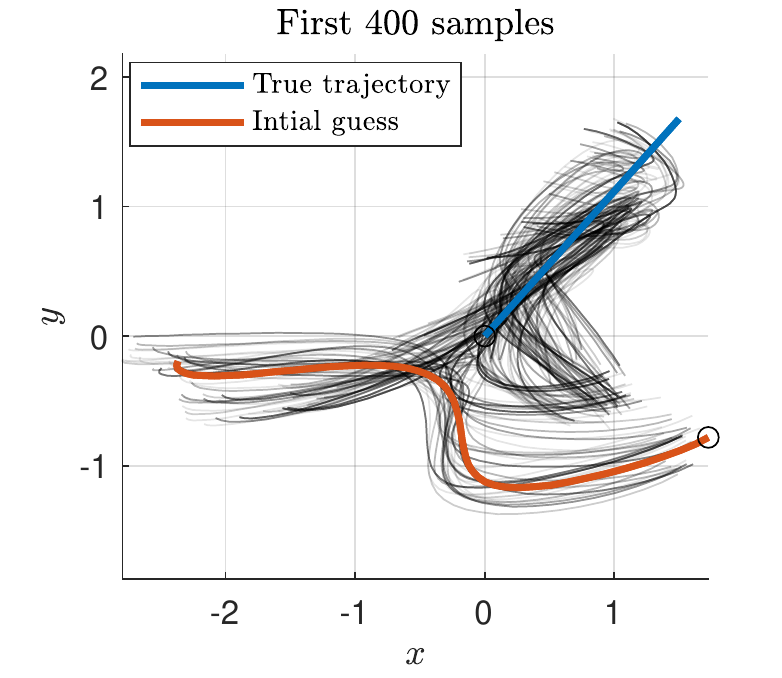}
$\ell = 4$, $\kappa = 1$
\end{minipage}
\begin{minipage}[c]{0.32\textwidth}
\centering
    \includegraphics[width=\linewidth]{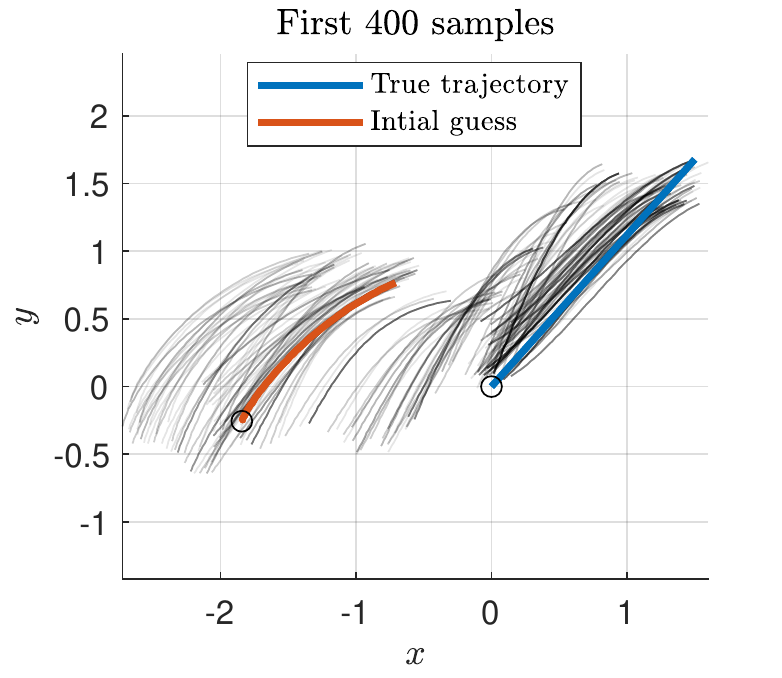}
$\ell = 15$, $\kappa = 1$
\end{minipage}
\caption{First 400 samples from one initial guess.}
\label{fig:case1_warmup}
\begin{minipage}[c]{0.32\textwidth}
\centering
    \includegraphics[width=\linewidth]{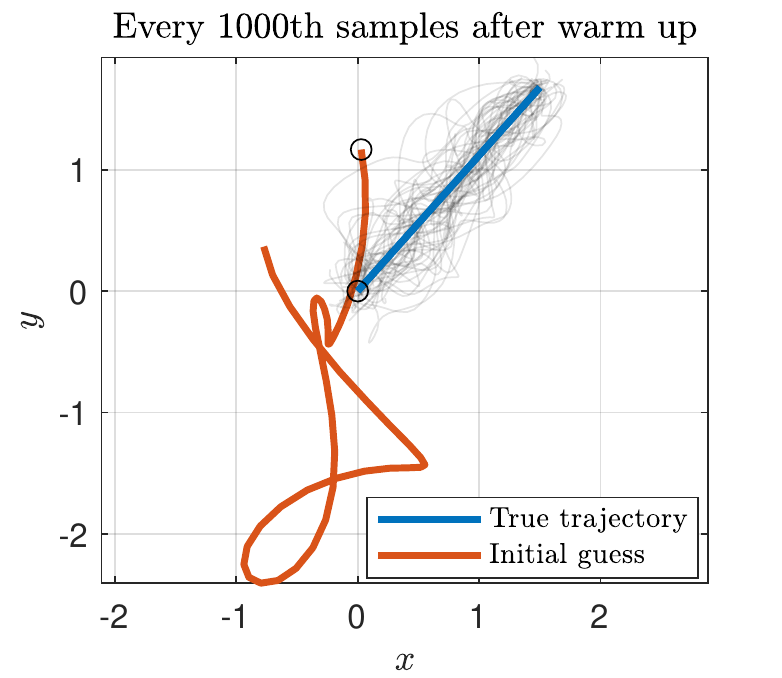}
$\ell = 2$, $\kappa = 1$
\end{minipage}
\begin{minipage}[c]{0.32\textwidth}
\centering
    \includegraphics[width=\linewidth]{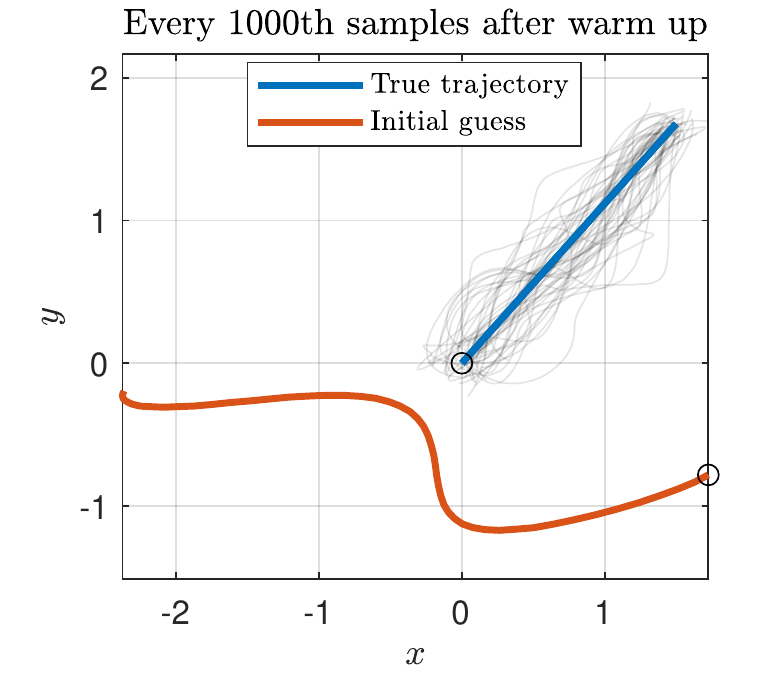}
$\ell = 4$, $\kappa = 1$
\end{minipage}
\begin{minipage}[c]{0.32\textwidth}
\centering
    \includegraphics[width=\linewidth]{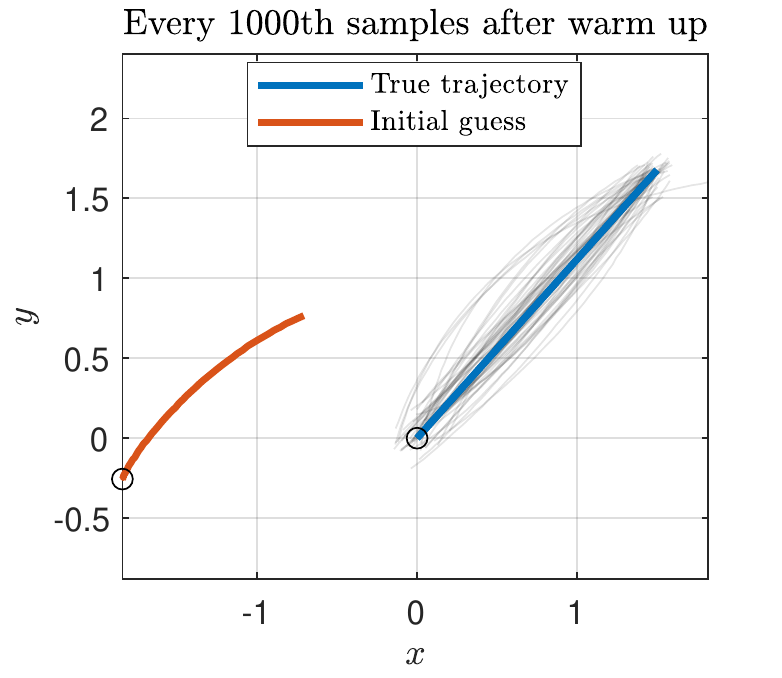}
$\ell = 15$, $\kappa = 1$
\end{minipage}
\caption{Every 1000th sample after warm-up from one initial guess.}
    \label{fig:case1_1000}
\end{figure}

\begin{figure}[h!]
\begin{minipage}[c]{0.32\textwidth}
\centering
    \includegraphics[width=\linewidth]{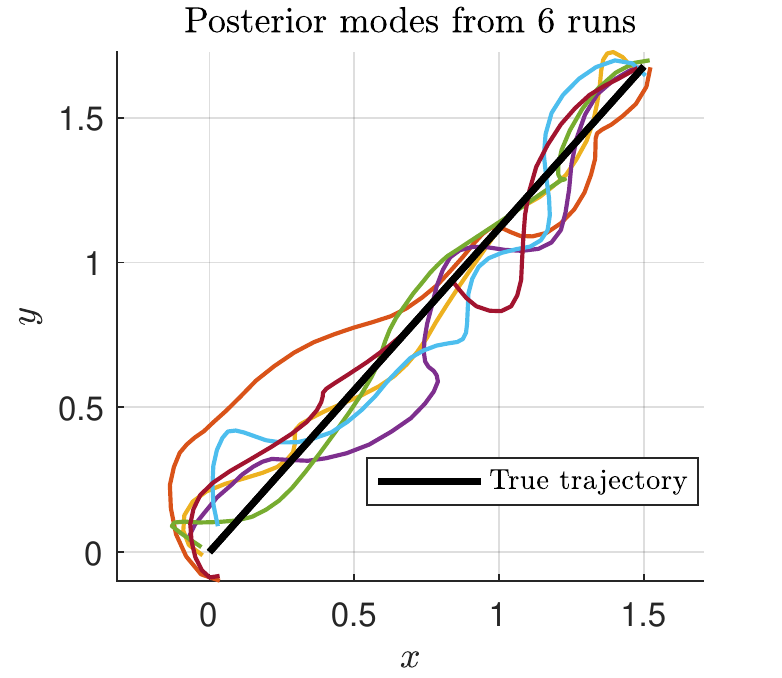}
$\ell = 2$, $\kappa =1$
\end{minipage}
\begin{minipage}[c]{0.32\textwidth}
\centering
    \includegraphics[width=\linewidth]{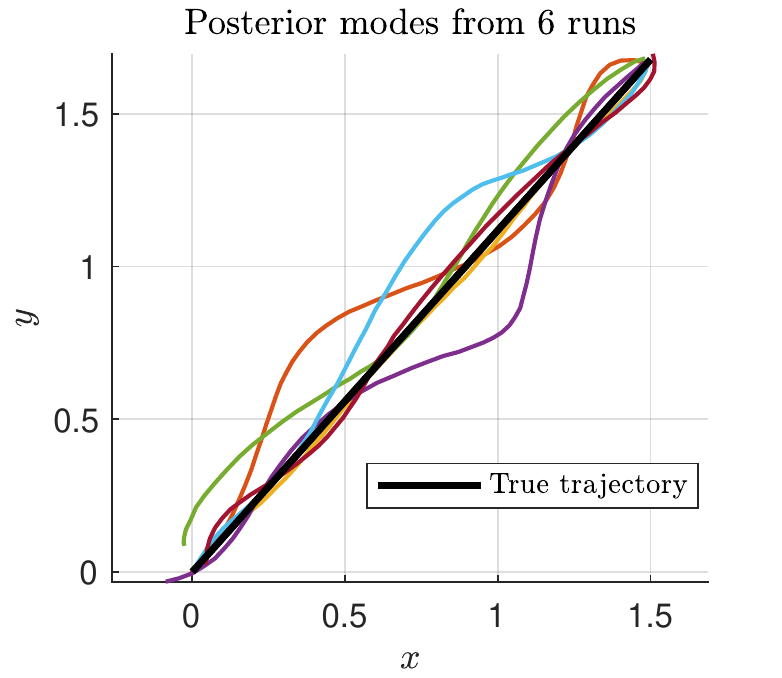}
$\ell = 4$, $\kappa = 1$
\end{minipage}
\begin{minipage}[c]{0.32\textwidth}
\centering
    \includegraphics[width=\linewidth]{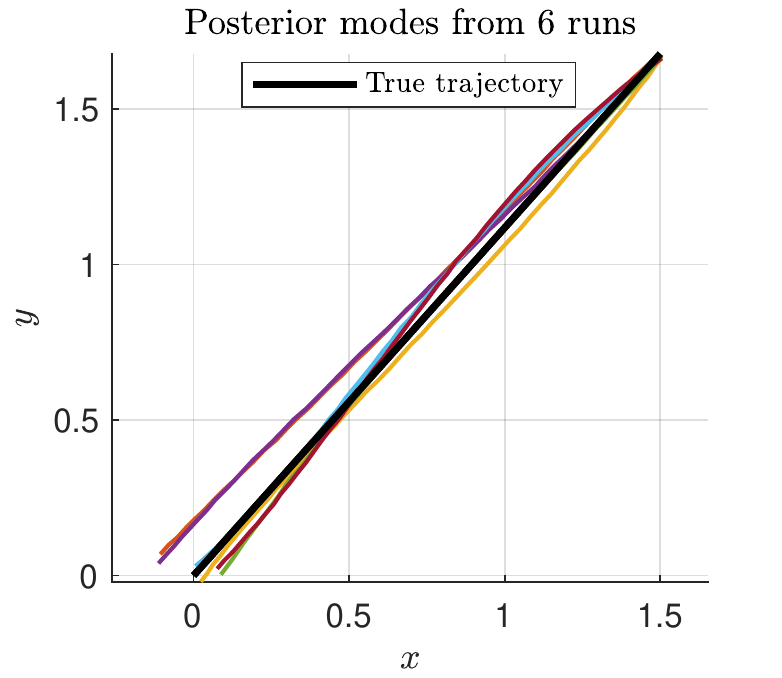}
$\ell = 15$, $\kappa = 1$
\end{minipage}
\caption{Posterior modes of the trajectory from six different initial guesses.}
\label{fig:case1_tra_mode}
\begin{minipage}[c]{0.32\textwidth}
\centering
    \includegraphics[width=\linewidth]{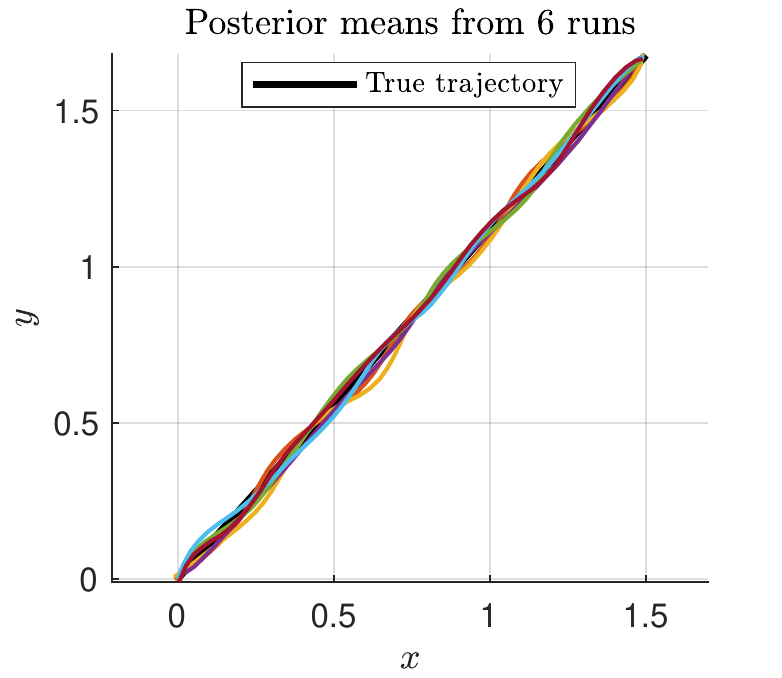}
$\ell = 2$, $\kappa =1$
\end{minipage}
\begin{minipage}[c]{0.32\textwidth}
\centering
    \includegraphics[width=\linewidth]{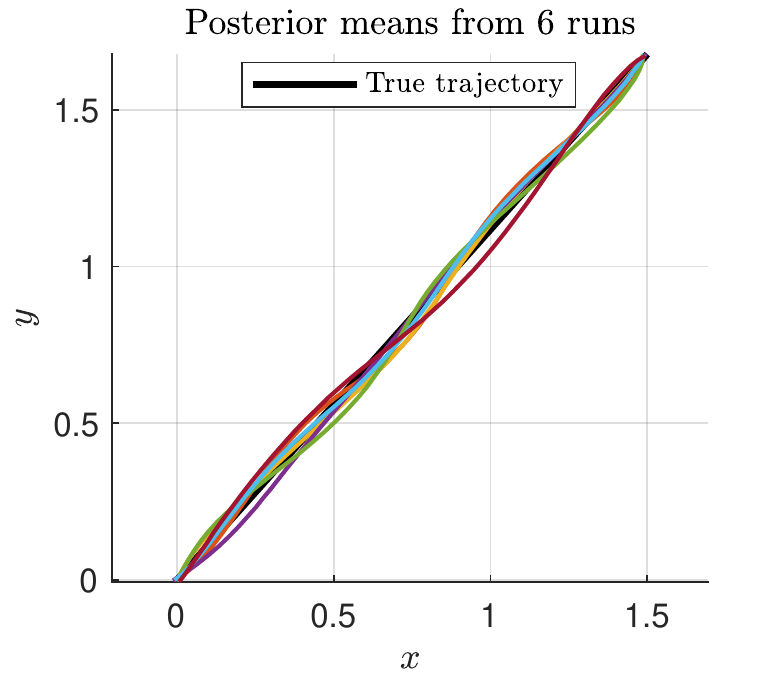}
$\ell = 4$, $\kappa = 1$
\end{minipage}
\begin{minipage}[c]{0.32\textwidth}
\centering
    \includegraphics[width=\linewidth]{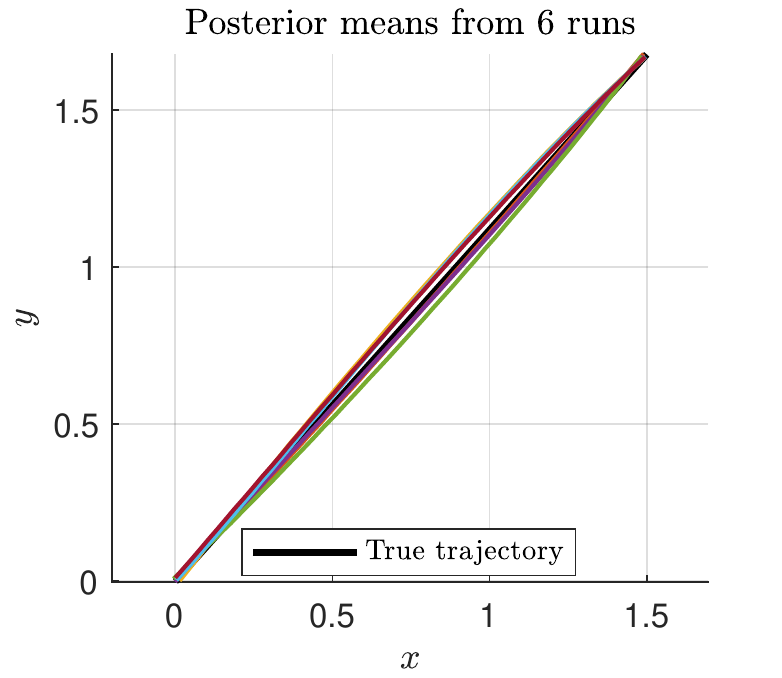}
$\ell = 15$, $\kappa = 1$
\end{minipage}
\caption{Posterior means of the trajectory from six different initial guesses.}\label{fig:case1_tra_mean}
\begin{minipage}[c]{0.32\textwidth}
\centering
    \includegraphics[width=\linewidth]{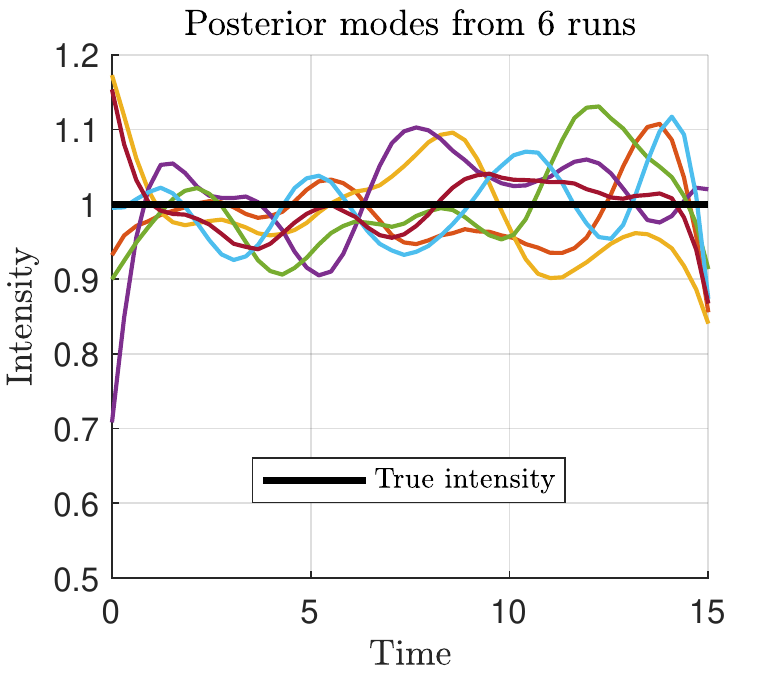}
$\ell = 2$, $\kappa =1$
\end{minipage}
\begin{minipage}[c]{0.32\textwidth}
\centering
    \includegraphics[width=\linewidth]{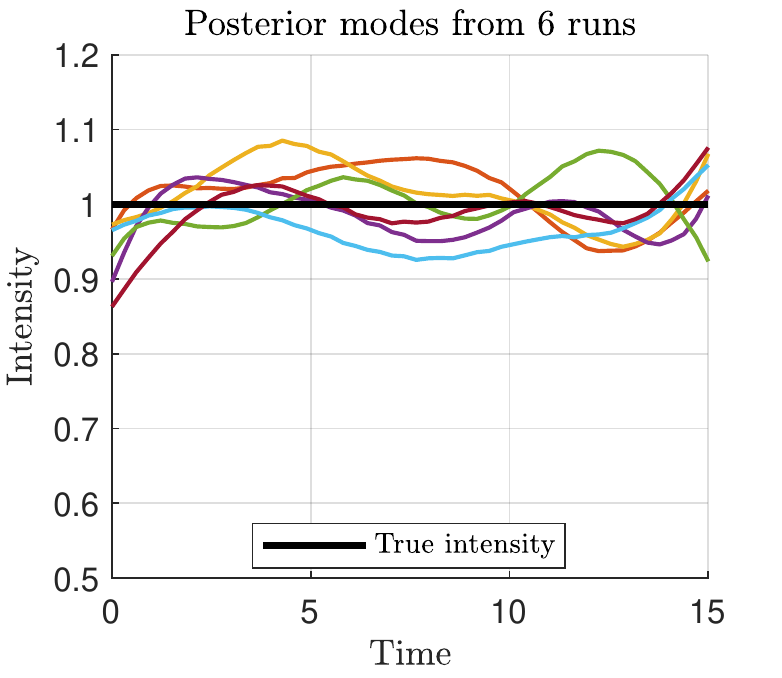}
$\ell = 4$, $\kappa = 1$
\end{minipage}
\begin{minipage}[c]{0.32\textwidth}
\centering
    \includegraphics[width=\linewidth]{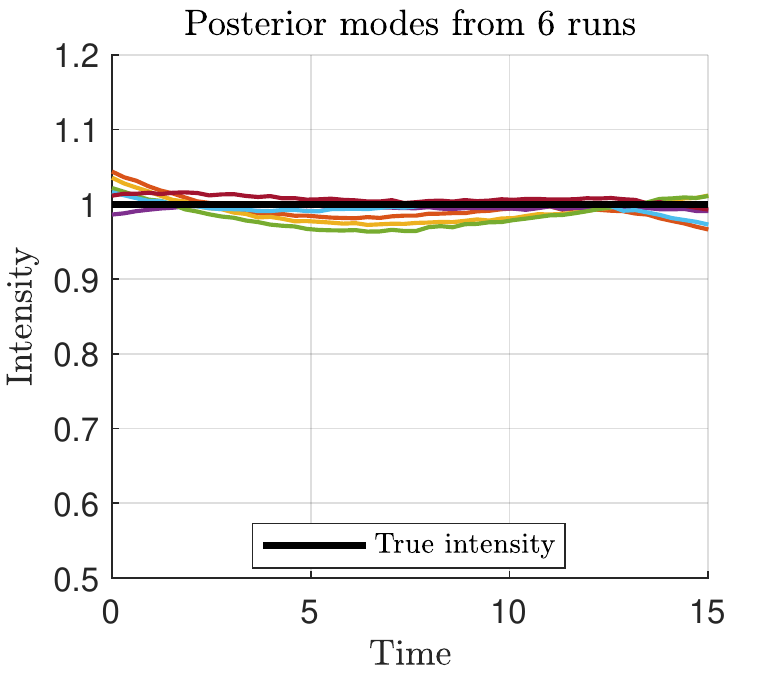}
$\ell = 15$, $\kappa = 1$
\end{minipage}
\caption{Posterior modes of the intensity from six different initial guesses.}\label{fig:case1_int_mode}
\begin{minipage}[c]{0.32\textwidth}
\centering
    \includegraphics[width=\linewidth]{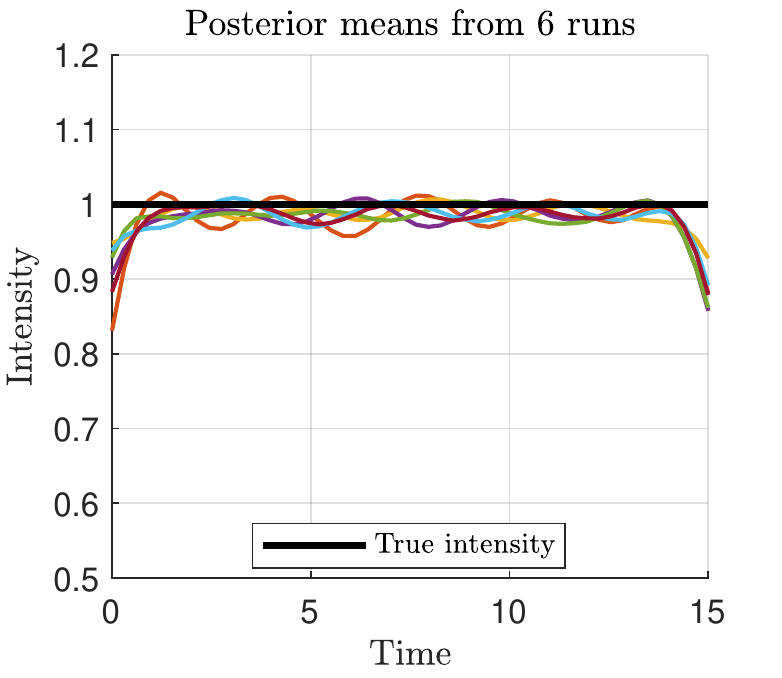}
$\ell = 2$, $\kappa =1$
\end{minipage}
\begin{minipage}[c]{0.32\textwidth}
\centering
    \includegraphics[width=\linewidth]{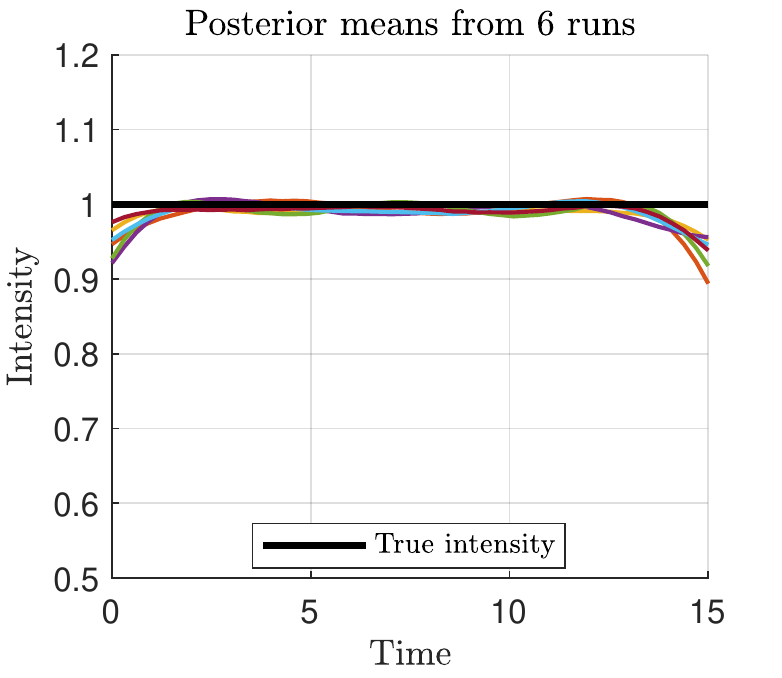}
$\ell = 4$, $\kappa = 1$
\end{minipage}
\begin{minipage}[c]{0.32\textwidth}
\centering
    \includegraphics[width=\linewidth]{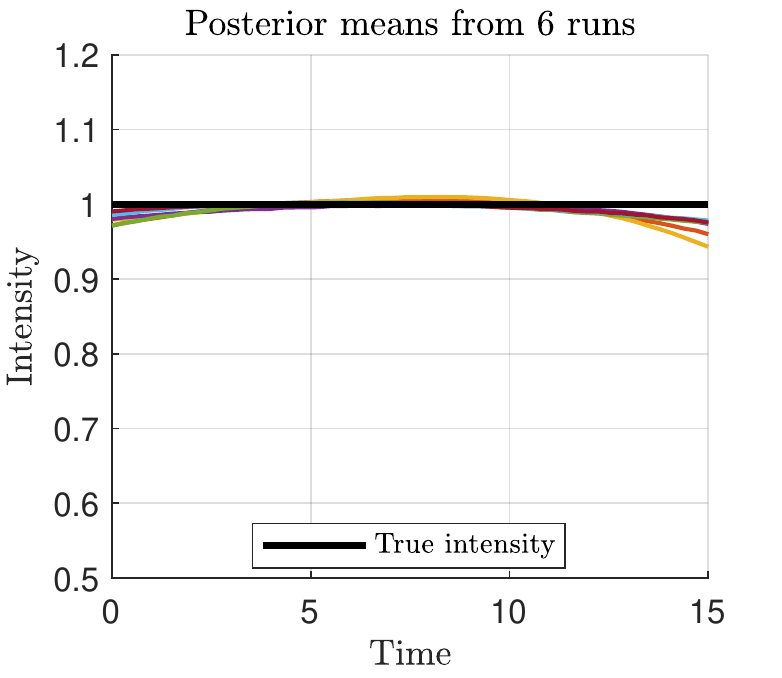}
$\ell = 15$, $\kappa = 1$
\end{minipage}
\caption{Posterior means of the intensity from six different initial guesses.}\label{fig:case1_int_mean}
\end{figure}

Straight-line trajectory in the $xy$-plane, $p(t)=(v_xt,v_yt,0)$, $t\in[0,T_0]$, with constant speed $v=\sqrt{v_x^2+v_y^2}=0.15$; constant intensity during the emission time, 
\[
q(t)=\begin{cases}1,\quad&t\in[0,T_0],\\0\quad&\text{otherwise.}\end{cases}
\]
The trajectory and intensity were to be estimated using 424 sensors uniformly distributed on a hemisphere. The measurements of the wave field are shown in Figure~\ref{fig:case1_measurements}.
As part of our study of Case 1, we investigated the effect of the hyperparameter $\ell$ in the covariance function on the numerical results. Figure~\ref{fig:case1_initial} shows the initial guesses for the source trajectory corresponding to three different values of $\ell$ and, as expected, larger values resulted in slower-turning trajectories
Figure~\ref{fig:case1_warmup} shows, for each of the three sets of hyperparameter values, the first 400 samples following one of the initial guesses. Clearly only the best-adapted choice, with the relatively large correlation length, results in apparent convergence towards the true source trajectory. As shown in Figure~\ref{fig:case1_1000}, the subsequent sampling improves the situation for the two other choices of hyperparameter values, but the best-adapted choice still seems to show better convergence to the true source trajectory. 
Figures~\ref{fig:case1_tra_mode},~\ref{fig:case1_tra_mean},~\ref{fig:case1_int_mode}, and~\ref{fig:case1_int_mean} show the modes (the trajectories with the highest posterior density) and means of the posterior trajectories and intensities from the six initial guesses. In all three cases, the posterior means seem to be better predictors than the posterior modes. In particular, the posterior means give a significant improvement over the modes when suboptimal hyperparameters are chosen.
%
%
%
%
Table~\ref{tab:case1_errors} shows the reconstruction errors, defined in Eqs. \eqref{eq:errors} for the three sets of values of the hyperparameters.
\begin{table}[h]
\begin{center}
  \begin{tabular}{ l | c | c | c | }
    & $\ell=2, \kappa=1$ & $\ell=4, \kappa=1$ & $\ell=15, \kappa=1$ \\ \hline
wavefield error & $9.5\cdot 10^{-4} $ & $7.8\cdot 10^{-4} $ & $6.0\cdot 10^{-4} $ \\ \hline
Average trajectory error & $5.1\cdot 10^{-2} $ & $4.5\cdot 10^{-2} $ & $4.1\cdot 10^{-2} $ \\ \hline
Average intensity error & $3.0 \cdot 10^{-2} $ & $1.9 \cdot 10^{-2} $ & $1.2\cdot 10^{-2} $ \\ \hline
  \end{tabular}
  \caption{Field, trajectory and intensity error for the average (over 6 starting guesses) posterior means for three sets of hyperparameter values.}\label{tab:case1_errors}
\end{center}
\end{table}

\FloatBarrier 
\subsection{Case 2} 
\begin{figure}[ht]
\centering
\begin{minipage}[c]{0.48\textwidth}
\centering
    \includegraphics[width=\linewidth]{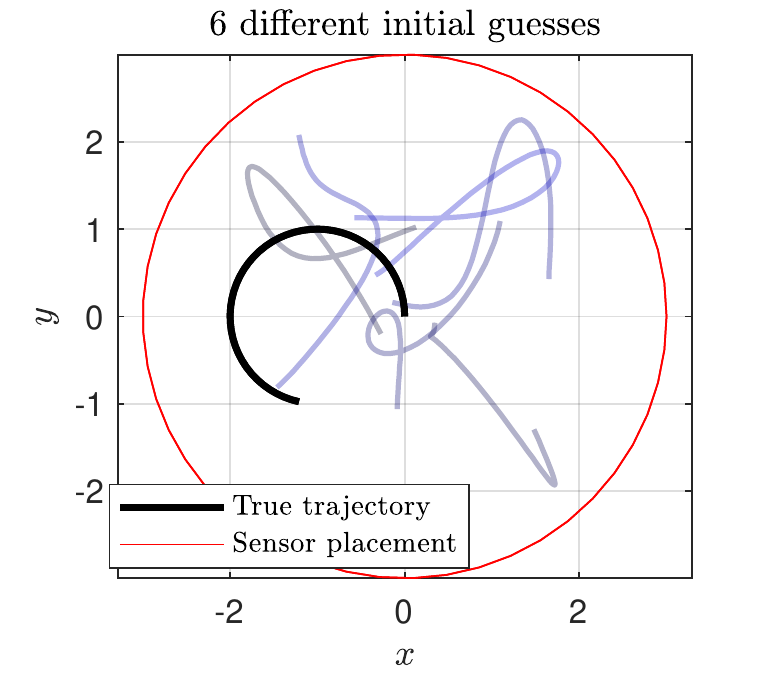}
    \end{minipage}
    \begin{minipage}[c]{0.48\textwidth}
\centering
        \includegraphics[width=\linewidth]{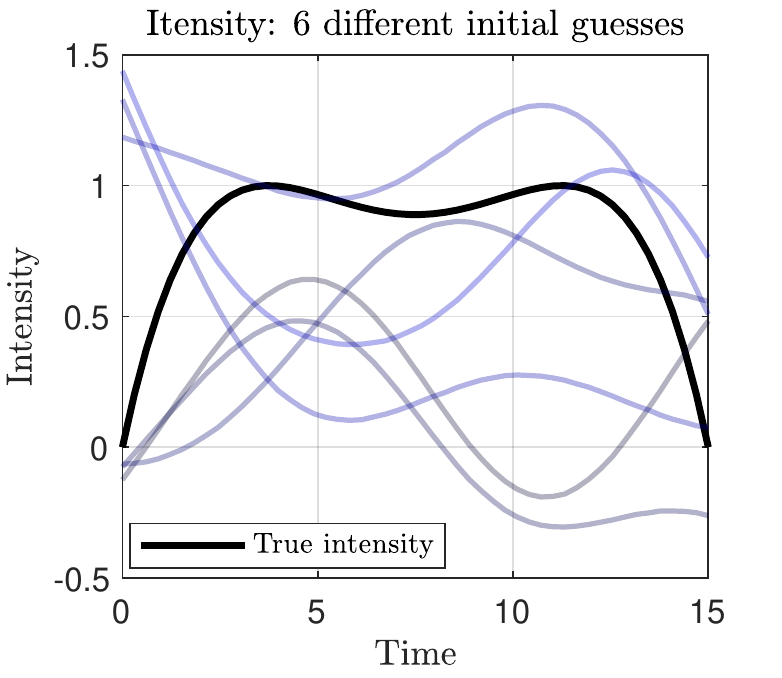}
\end{minipage}
\caption{Initial samples from the GP prior for the trajectory (left) and intensity (right).}
    \label{fig:circle_initial_guess} 
\end{figure}

\begin{figure}[ht]
\centering
\begin{minipage}[c]{0.48\textwidth}
\centering
    \includegraphics[width=\linewidth]{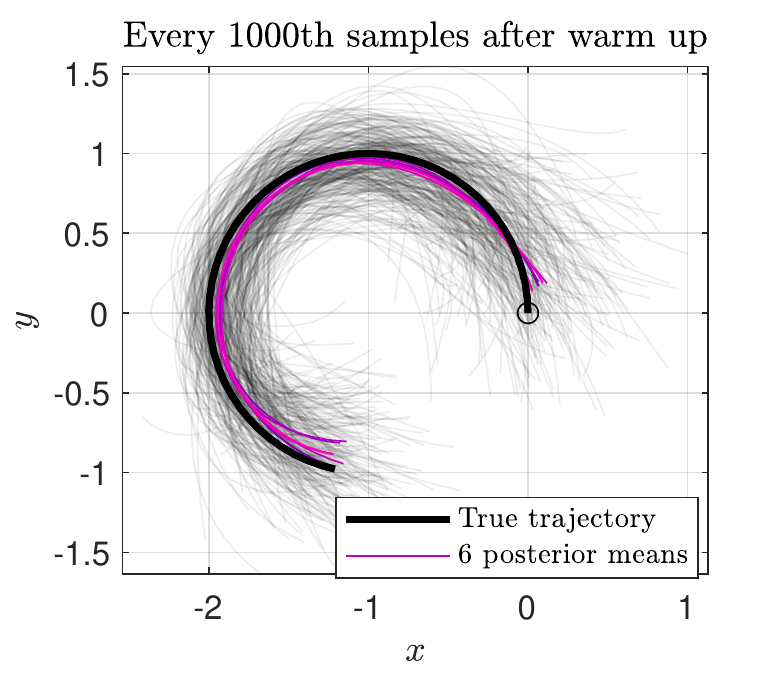}
    $(a)$
    \end{minipage}
    \begin{minipage}[c]{0.48\textwidth}
\centering
        \includegraphics[width=\linewidth]{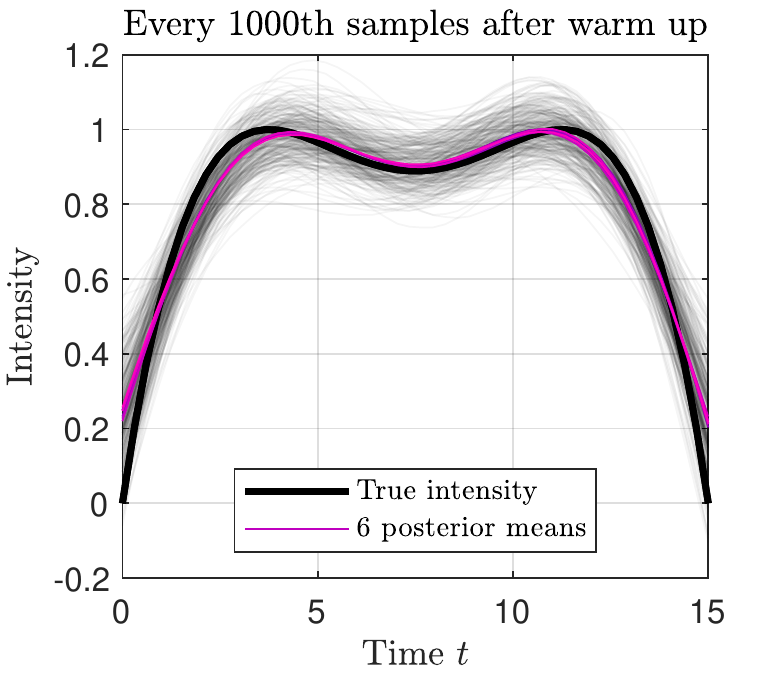}
        $(b)$
\end{minipage}
\caption{Every 1000th sample after warm of from all 6 chains, and their corresponding posterior mean for $(a)$ the trajectory and $(b)$ the intensity.}
\label{fig:circle_every_1000th} 
\end{figure}

We consider a circular arc trajectory in the $xy$-plane parametrized by 
\begin{equation}
    p(t)=(\cos(0.3t)-1,\sin(0.3t),0), t\in[0,T_0]
\end{equation}
and with non-constant polynomial intensity
\begin{equation}
    q(t)=\begin{cases}
    -28.44(t/T_0)^4+56.89(t/T_0)^3-39.11(t/T_0)^2+10.67(t/T_0), &t\in[0,T_0] \\
    0, & \text{otherwise}
    \end{cases}.
\end{equation}
We use the hyperparameter values $\kappa=1$, $\ell=5$ with 6 resulting samples from the GP prior shown in Fig.~\ref{fig:circle_initial_guess}. Figure~\ref{fig:circle_every_1000th} shows every 1000th sample from six independent MCMC chains together with the corresponding posterior means. 

\begin{figure}
\begin{minipage}[c]{0.45\textwidth}
\centering
    \includegraphics[width=\linewidth]{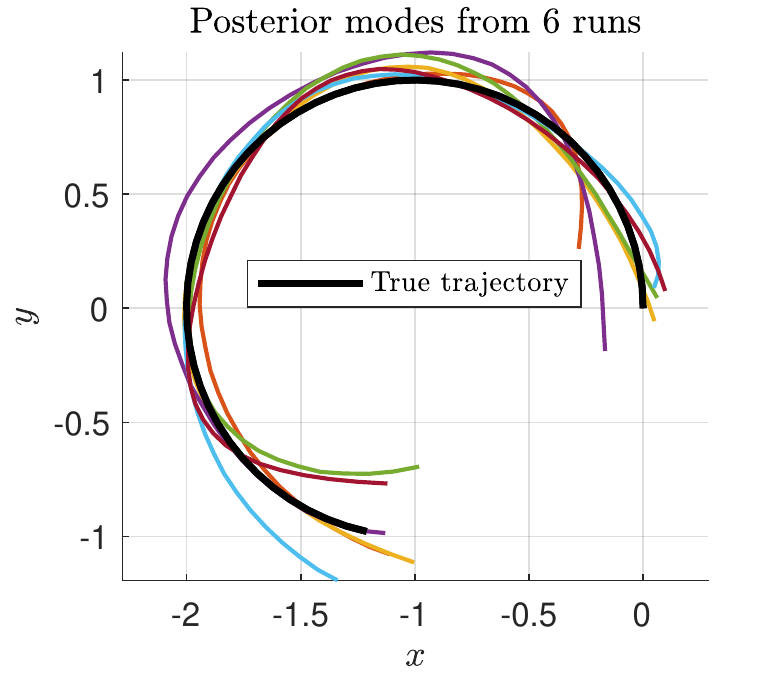}
    \includegraphics[width=\linewidth]{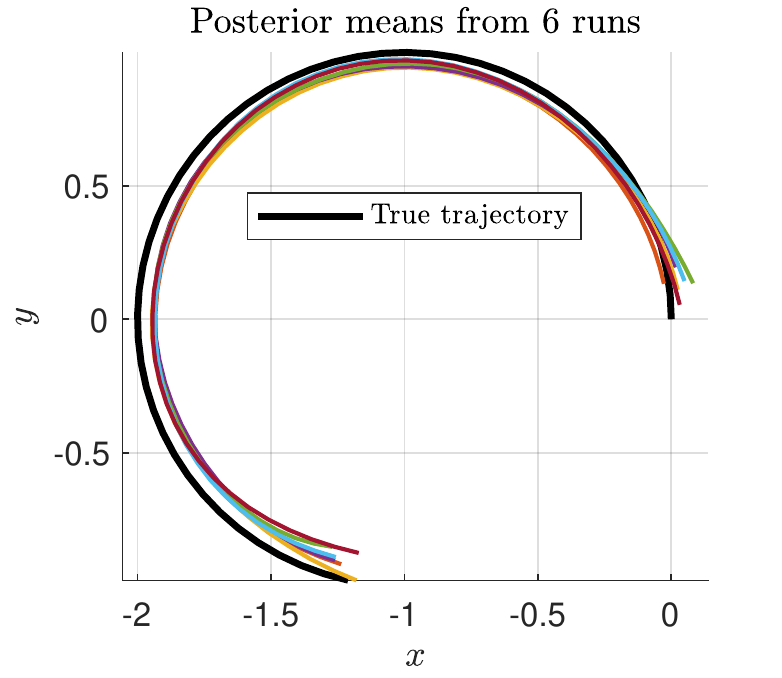}
        \includegraphics[width=\linewidth]{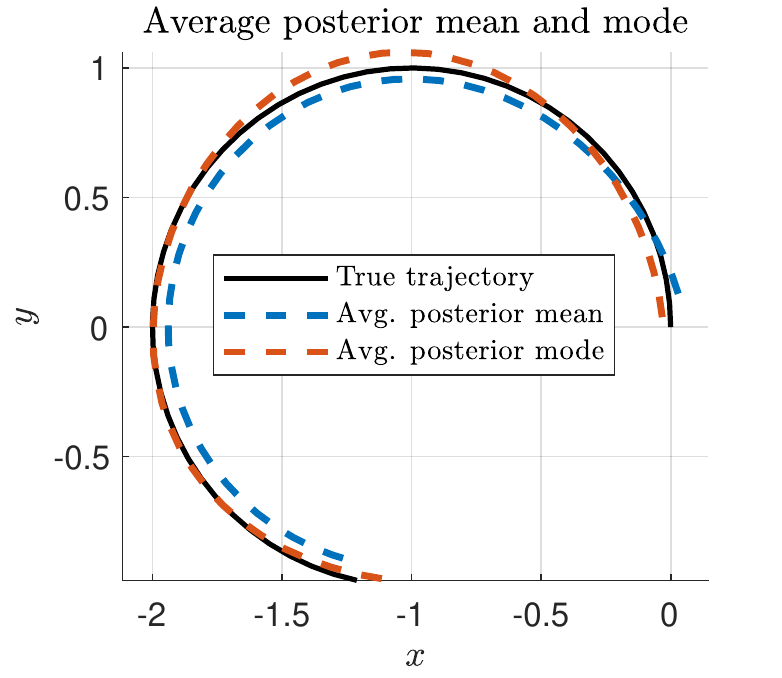}
\end{minipage}
\begin{minipage}[c]{0.45\textwidth}
\centering
    \includegraphics[width=\linewidth]{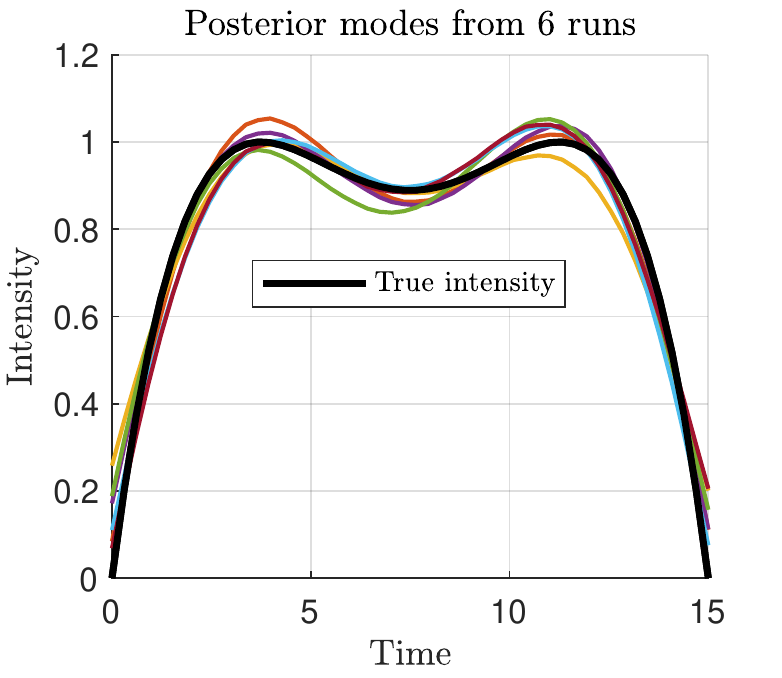}
    \includegraphics[width=\linewidth]{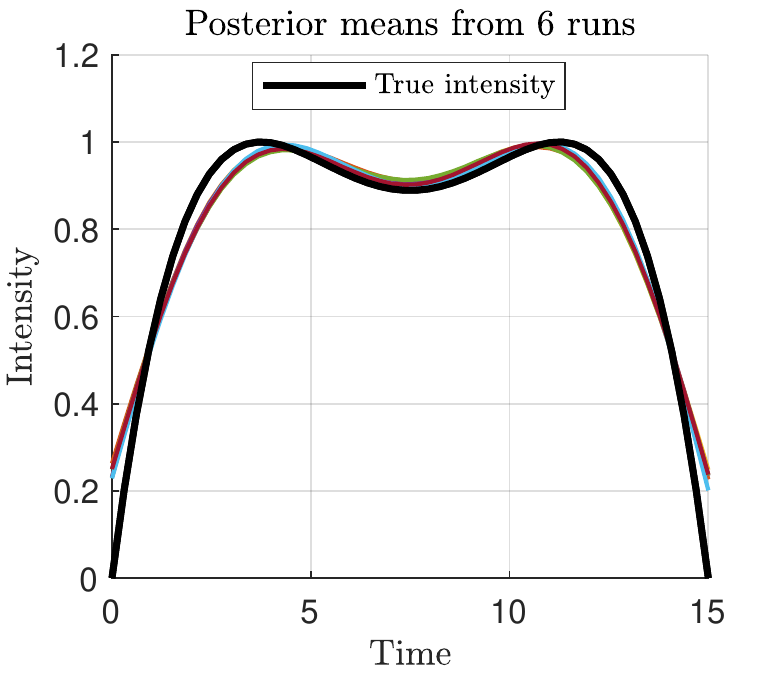}
        \includegraphics[width=\linewidth]{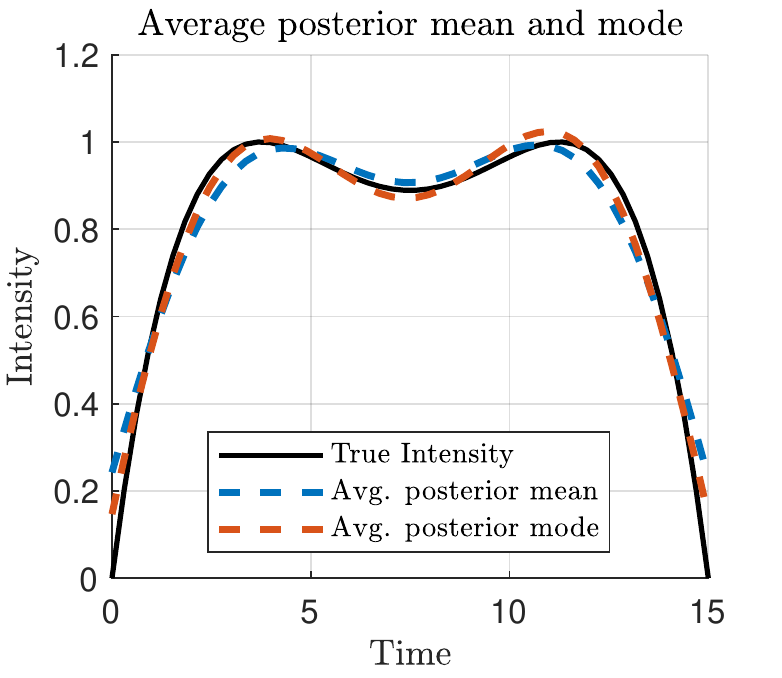}
\end{minipage}
\caption{Reconstruction with polynomial intensity.}
\label{fig:circle_recon_pol_intensity}
\end{figure}

\begin{figure}
\begin{minipage}[c]{0.32\textwidth}
\centering
    \includegraphics[width=\linewidth]{fig/circle/30_b.pdf}
\end{minipage}
\begin{minipage}[c]{0.32\textwidth}
\centering
    \includegraphics[width=\linewidth]{fig/circle/30_c.pdf}
\end{minipage}
\begin{minipage}[c]{0.32\textwidth}
\centering
    \includegraphics[width=\linewidth]{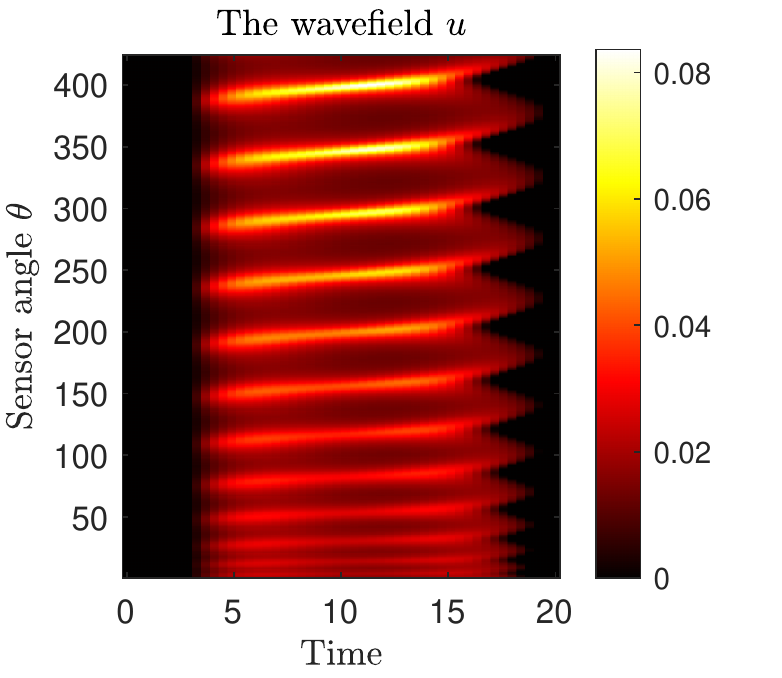}
\end{minipage}
\begin{center}
 No noise
\end{center}
\begin{minipage}[c]{0.32\textwidth}
\centering
    \includegraphics[width=\linewidth]{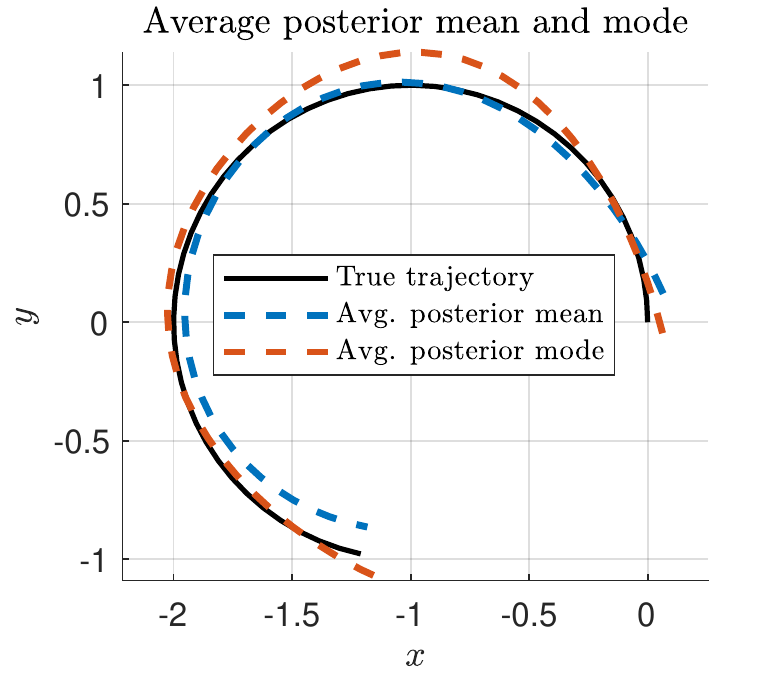}
\end{minipage}
\begin{minipage}[c]{0.32\textwidth}
\centering
    \includegraphics[width=\linewidth]{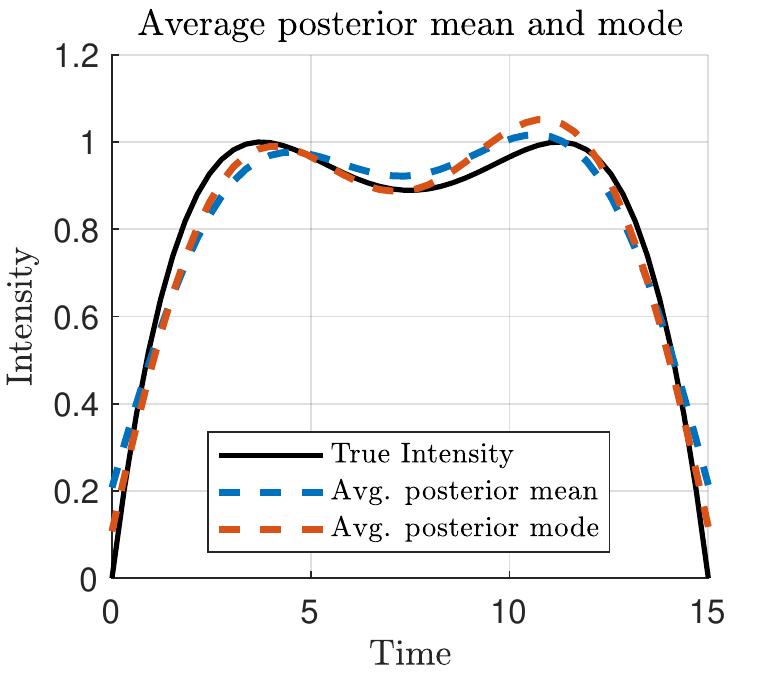}
\end{minipage}
\begin{minipage}[c]{0.32\textwidth}
\centering
    \includegraphics[width=\linewidth]{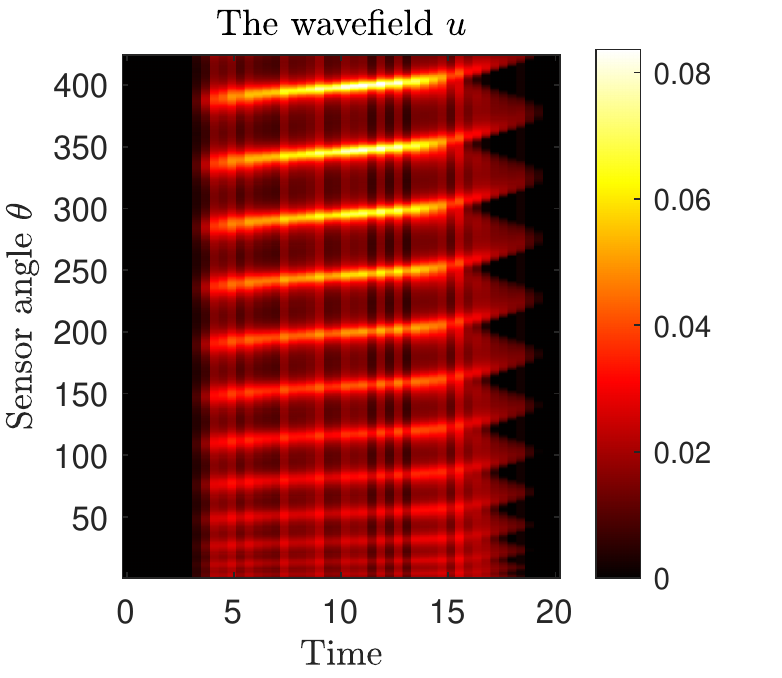}
\end{minipage}
\begin{center}
 $5\%$ noise
\end{center}
\begin{minipage}[c]{0.32\textwidth}
\centering
    \includegraphics[width=\linewidth]{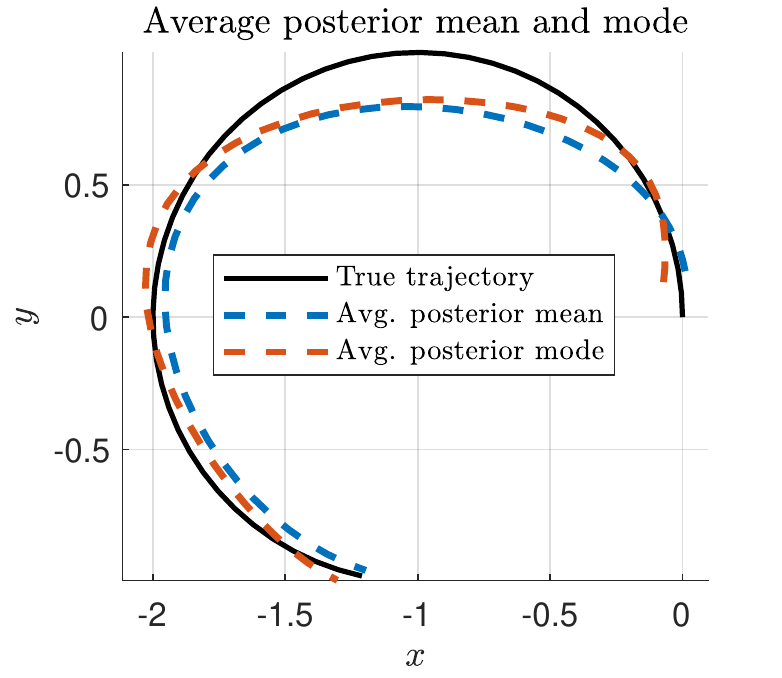}
\end{minipage}
\begin{minipage}[c]{0.32\textwidth}
\centering
    \includegraphics[width=\linewidth]{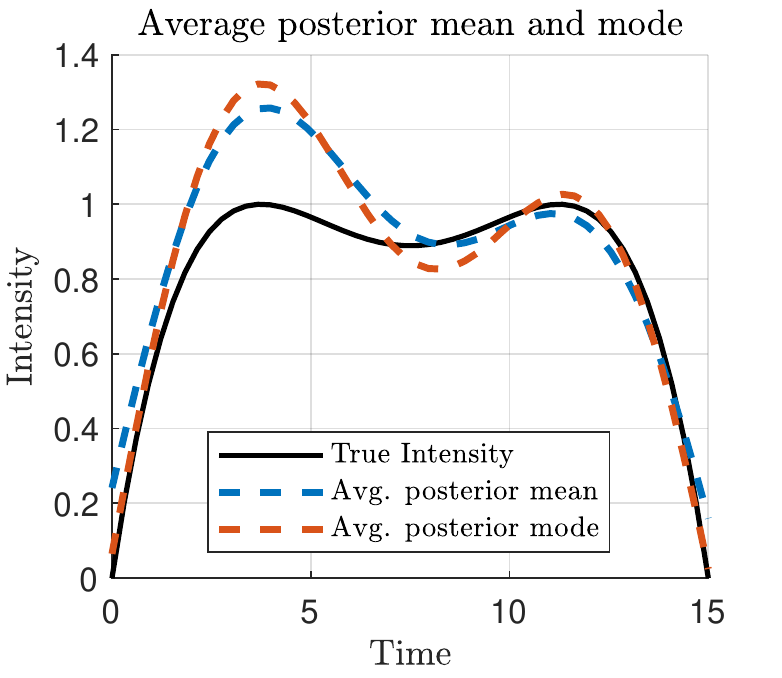}
\end{minipage}
\begin{minipage}[c]{0.32\textwidth}
\centering
    \includegraphics[width=\linewidth]{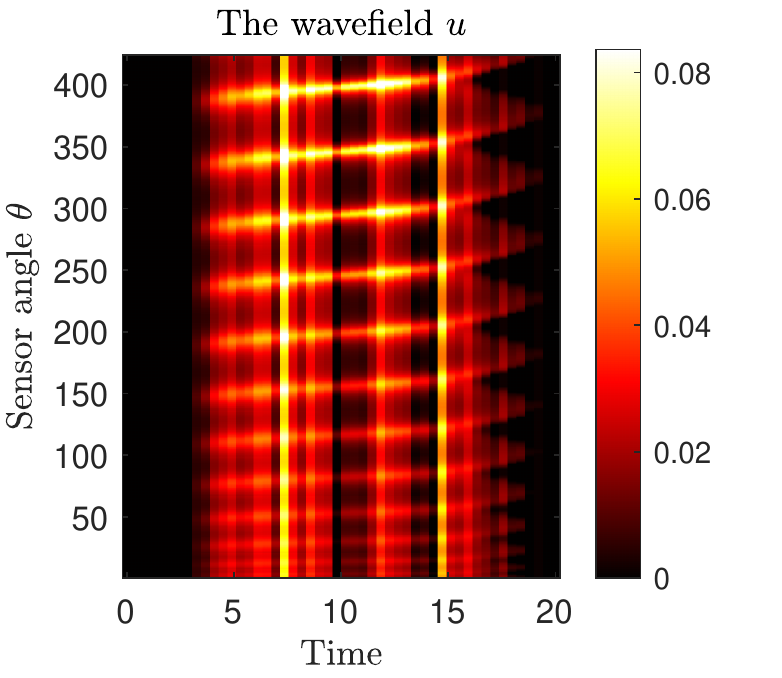}
\end{minipage}
\begin{center}
 $25\%$  noise
\end{center}
\begin{minipage}[c]{0.32\textwidth}
\centering
    \includegraphics[width=\linewidth]{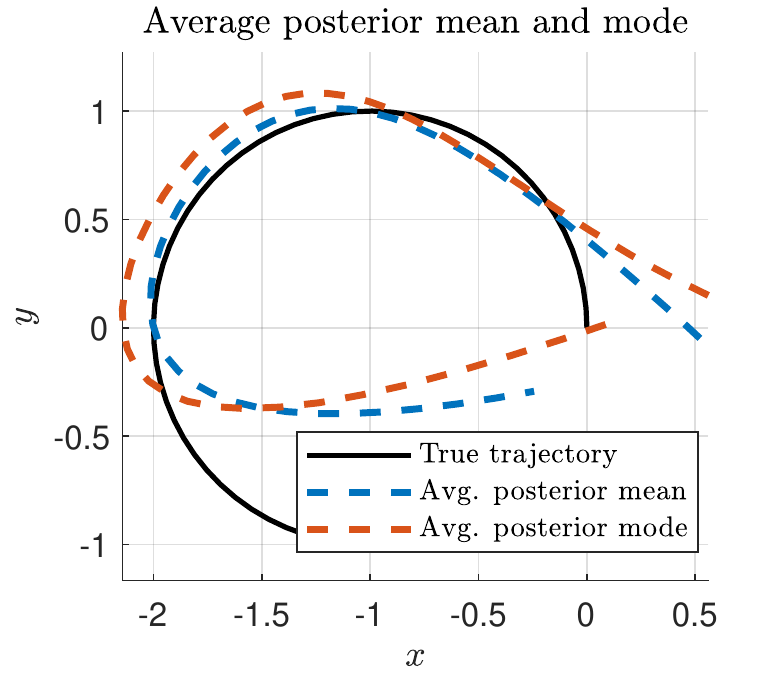}
\end{minipage}
\begin{minipage}[c]{0.32\textwidth}
\centering
    \includegraphics[width=\linewidth]{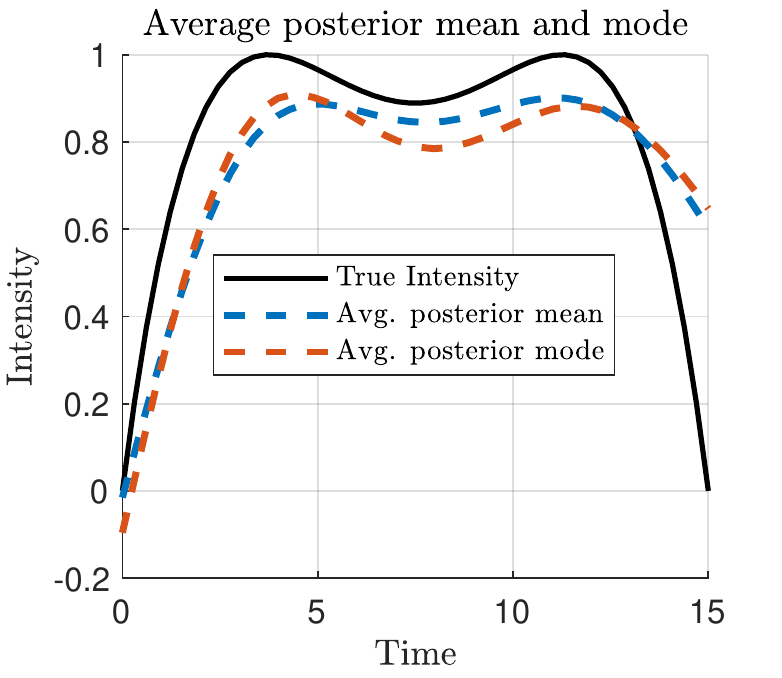}
\end{minipage}
\begin{minipage}[c]{0.32\textwidth}
\centering
    \includegraphics[width=\linewidth]{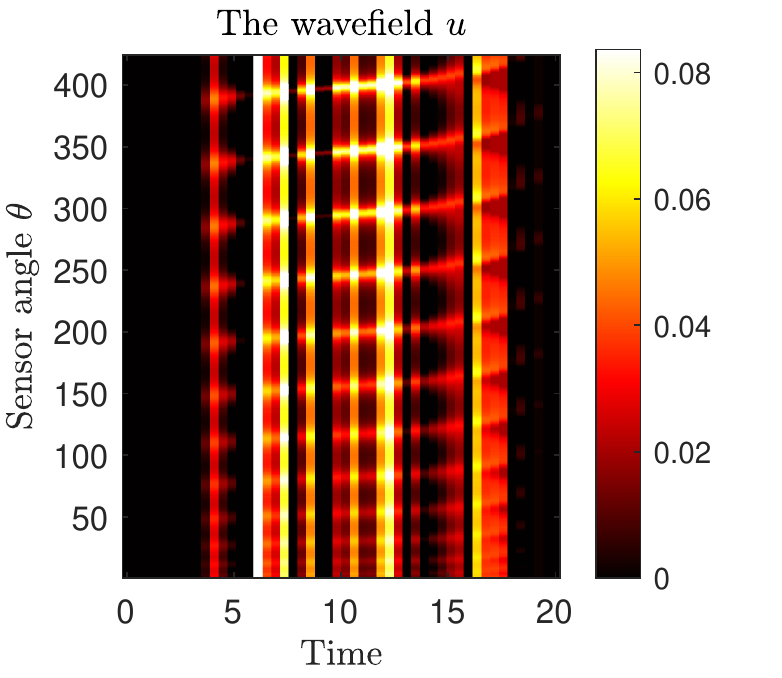}
\end{minipage}
\begin{center}
 $50\%$  noise
\end{center}
\caption{Reconstruction with polynomial intensity with no noise, $5\%$, $25\%$ and $50\%$ noise.}
\label{fig:circle_recon_pol_intensity_noise}
\end{figure}

A comparison between the posterior mean, posterior mode and averaged mean and mode over the six chains are shown in Fig. \ref{fig:circle_recon_pol_intensity}.
The effective sample size of $N_\mathrm{eff} = 656$ is sufficient to represent the posterior distribution. We again see larger variance in the posterior modes than the more conservative posterior means and a slight bias towards origo for the means, possibly due to a slight asymmetry in the posterior distribution. The averaged means and modes give the best reconstruction as quantified in Tbl. \ref{tab:circle_mean_error}. We note that the reconstruction is worst near the trajectory end points due to the low emission intensity.
\begin{table}
\begin{center}
  \begin{tabular}{ l | c | c |  }
& Average mean & Average mode\\ \hline
wavefield error  & $1.4\cdot 10^{-3} $ & $8.9\cdot 10^{-4} $ \\ \hline
trajectory error &$5.9\cdot 10^{-2} $  & $4.9\cdot 10^{-2} $ \\ \hline
intensity error & $6.5 \cdot 10^{-2} $ & $3.9\cdot 10^{-2} $  \\ \hline
  \end{tabular}
  \caption{Error for the average posterior mean and the average posterior mode.}
  \label{tab:circle_mean_error}
\end{center}
\end{table}
We now consider the effect of noisy measurements on the reconstruction. We incorporate noise similarly to \cite{ohe} where at each measurement time $t = t_l$ we set $u_\mathrm{meas}(x,t_l) = u(x,t_l) + \epsilon_l$ where $u$ is the noise-free solution with additive time-dependent noise $\epsilon \sim \mathcal{N}(0,\sigma_l^2)$. The noise magnitude is defined relative to the field
\begin{equation}
    \sigma_l = \alpha \sqrt{\int_\Gamma |u(x,t_l)|^2 \mathrm{d}S} \approx \alpha \sqrt{\frac{A(\Gamma)}{N_\mathrm{sensors}}\sum_{i=1}^{N_\mathrm{sensors}} |u(x_i,t_l)|^2},
\end{equation}
where $A(\Gamma)$ is the area of the measurement surface. Note that this type of noise is purely time-dependent and is the same across the sensors, which is a harder problem than the case where the noise is also random between the sensors. The resulting reconstructions for $\alpha \in (0,0.05, 0.25, 0.5)$ are shown in Fig. \ref{fig:circle_recon_pol_intensity_noise}.
%
We obtain a reasonable reconstruction up to 25 \% noise with the intensity reconstruction being more affected than the trajectory. Our method seems to be highly resistant to noise due to using whole-trajectory samples compared to algebraic methods where the trajectory is reconstructed pointwise \cite{ohe}. However, direct comparison in difficult since \cite{ohe} considers three simultaneous sources for reconstruction.

\FloatBarrier
\subsection{Case 3}

\begin{figure}[ht]
\centering
\begin{minipage}[c]{0.48\textwidth}
\centering
    \includegraphics[width=\linewidth]{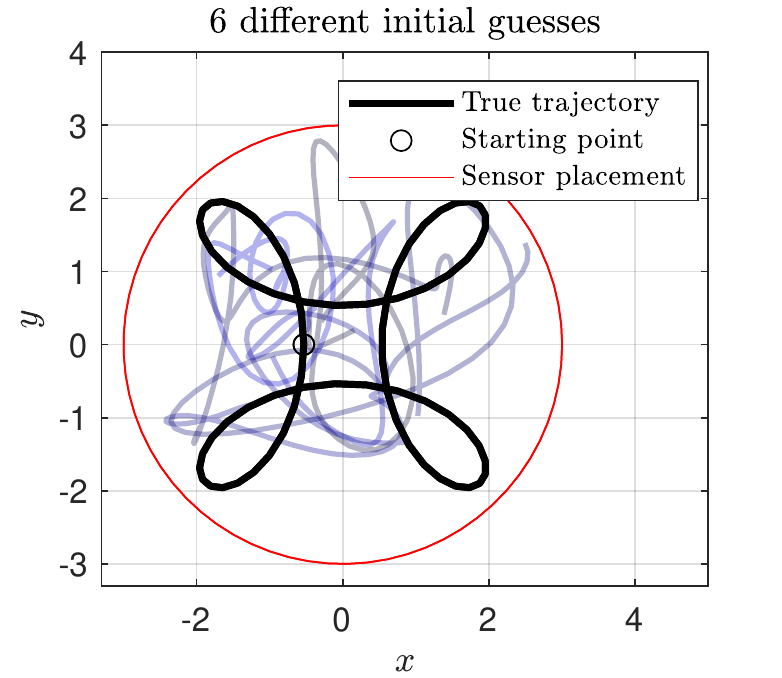}
    \end{minipage}
    \begin{minipage}[c]{0.48\textwidth}
\centering
        \includegraphics[width=\linewidth]{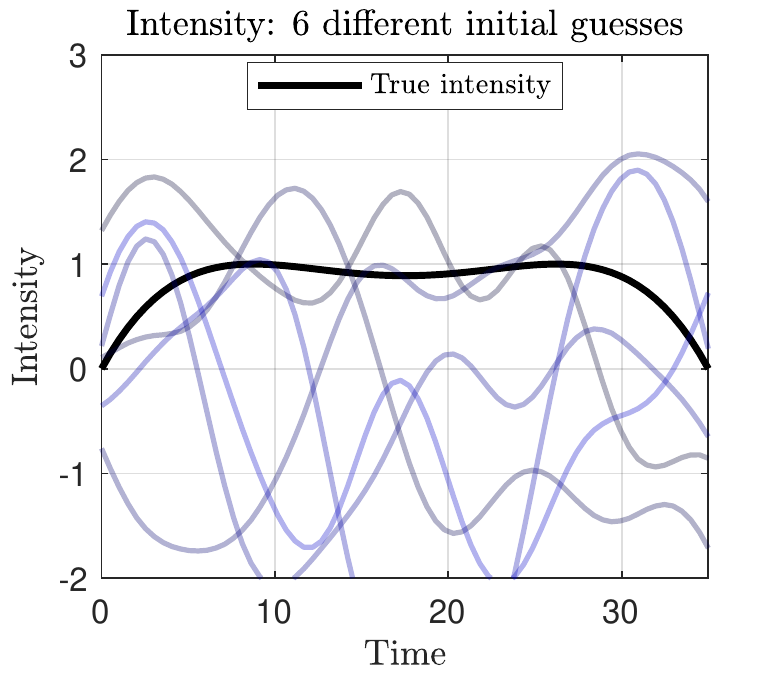}
\end{minipage}
\caption{The setup for the long trajectory. The source initially travels in the positive $y$-direction.}
    \label{fig:bow_set_up} 
\end{figure}

\begin{figure}[ht]
\centering
\begin{minipage}[c]{0.48\textwidth}
\centering
    \includegraphics[width=\linewidth]{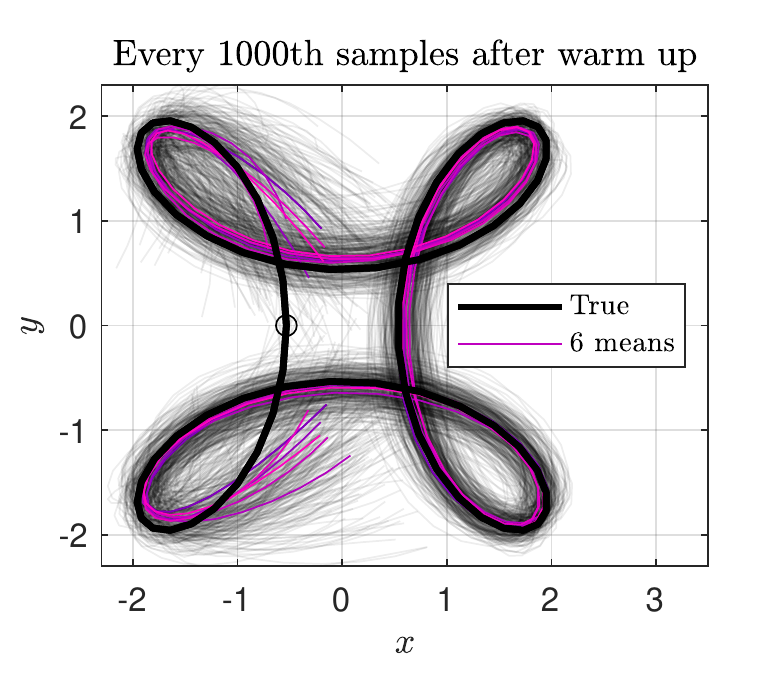}
    \end{minipage}
    \begin{minipage}[c]{0.48\textwidth}
\centering
        \includegraphics[width=\linewidth]{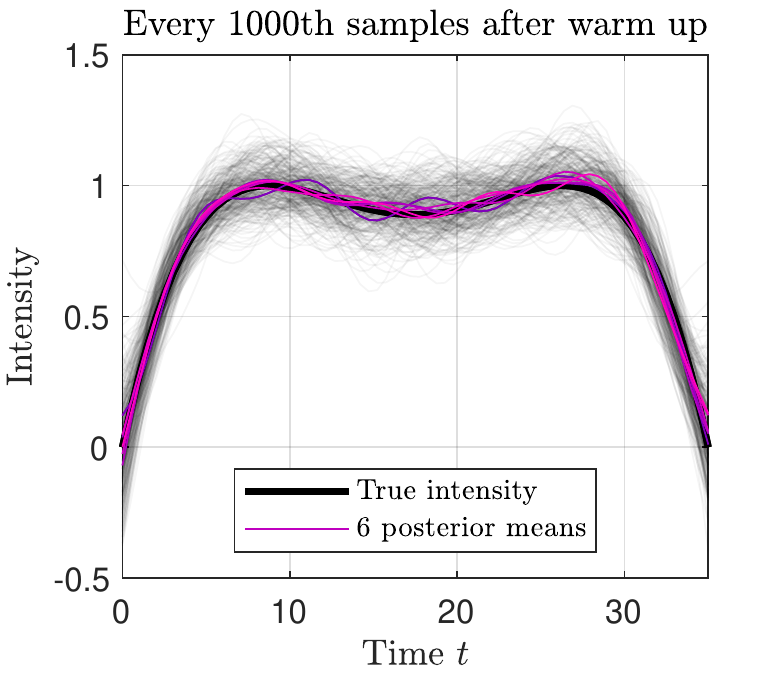}
\end{minipage}
\caption{Every 1000th sample from all 6 chains and their corresponding posterior mean.}
\label{fig:bow_every_1000th} 
\begin{minipage}[c]{0.45\textwidth}
\centering
    \includegraphics[width=\linewidth]{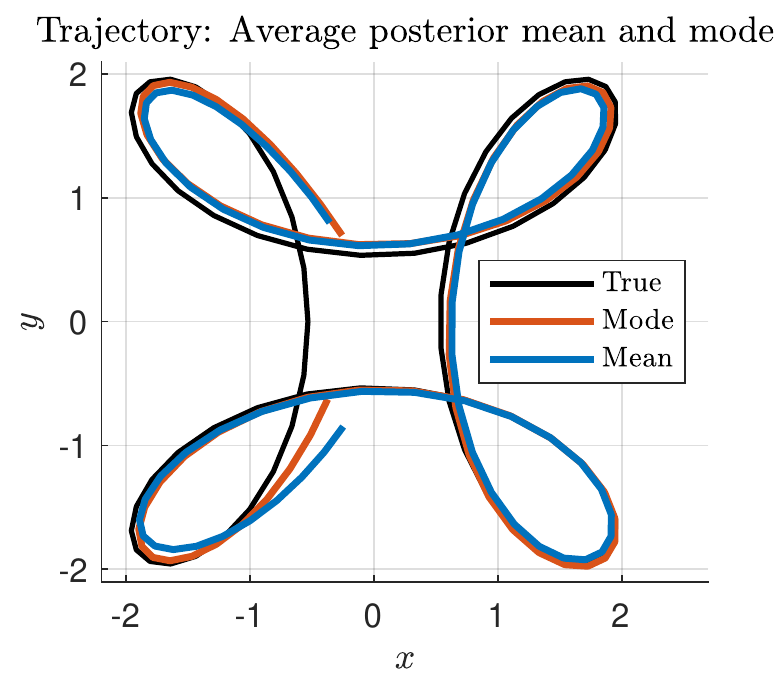}
\end{minipage}
\begin{minipage}[c]{0.45\textwidth}
\centering
    \includegraphics[width=\linewidth]{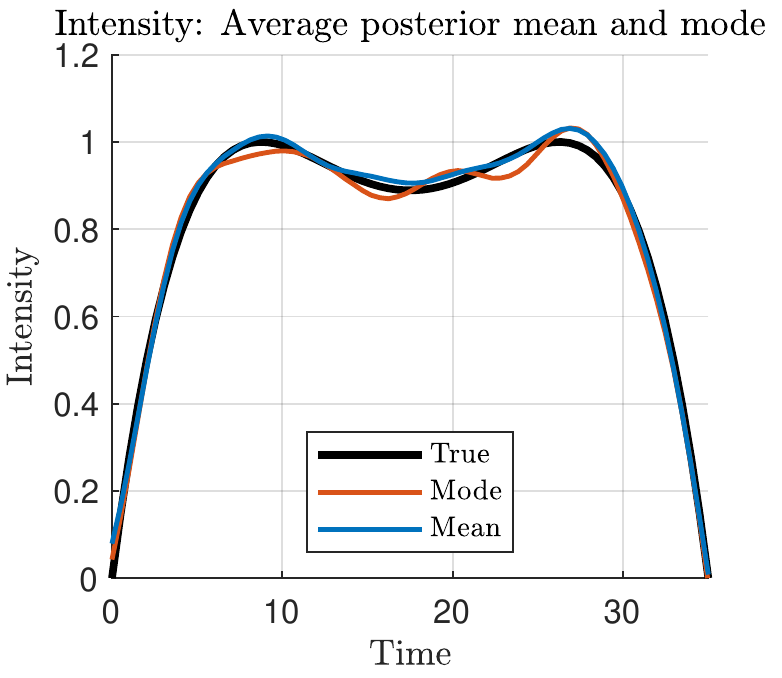}
\end{minipage}
\caption{The average posterior mean and mode of the trajectory and intensity.}
\label{fig:bow_average_mode_mean} 
\end{figure}

\begin{figure}[ht]
\centering
\begin{minipage}[c]{0.48\textwidth}
\centering
    \includegraphics[width=\linewidth]{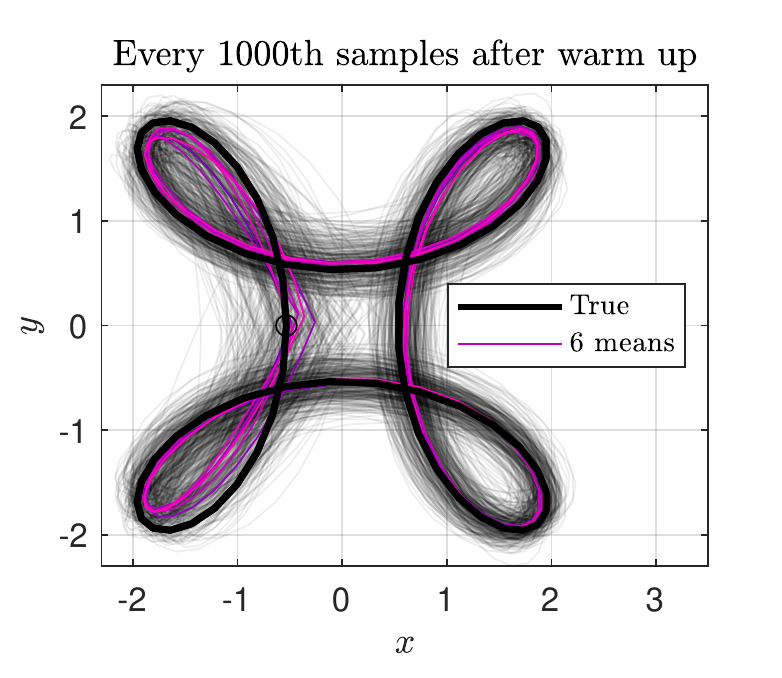}
    \end{minipage}
    \begin{minipage}[c]{0.48\textwidth}
\centering
        \includegraphics[width=\linewidth]{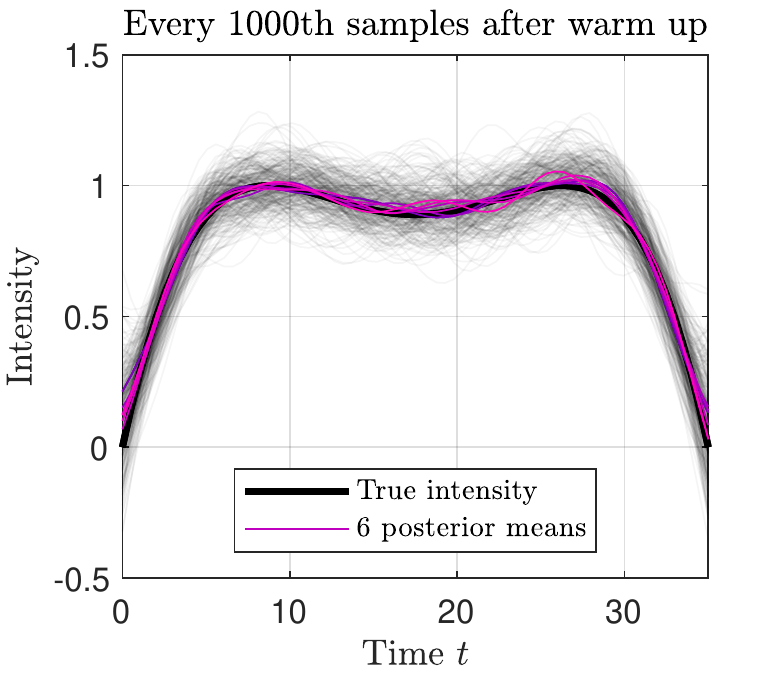}
\end{minipage}
\caption{Every 1000th sample from 6 chains and the corresponding posterior means with closed-curve conditioning.}
\label{fig:bow_every_1000th_condition} 
\begin{minipage}[c]{0.45\textwidth}
\centering
    \includegraphics[width=\linewidth]{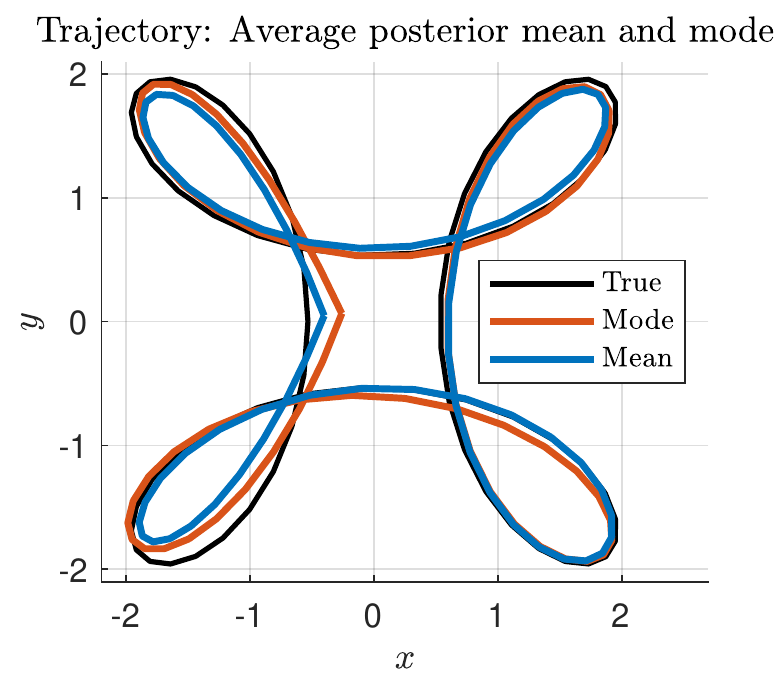}
\end{minipage}
\begin{minipage}[c]{0.45\textwidth}
\centering
    \includegraphics[width=\linewidth]{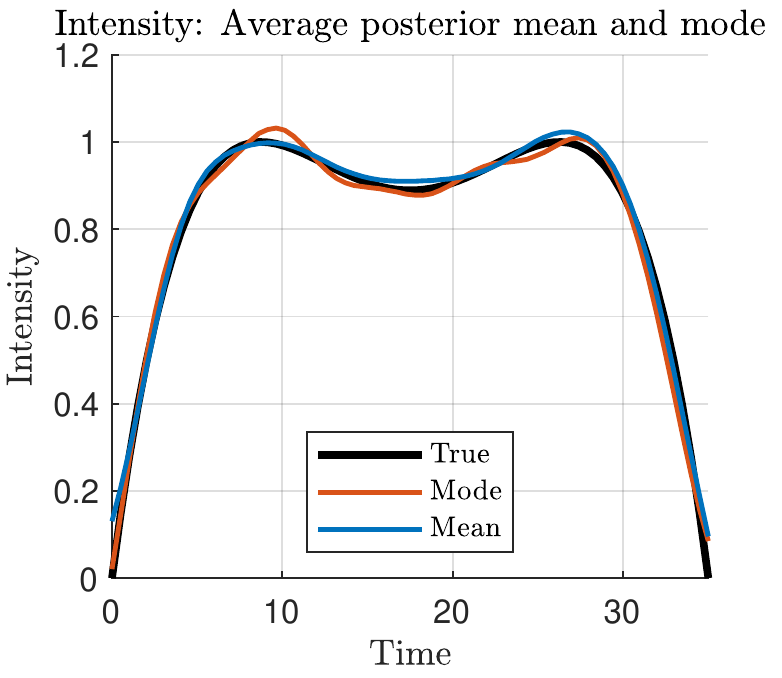}
\end{minipage}
\caption{The average posterior mean and mode of the trajectory and intensity with closed-curve conditioning.}
\label{fig:bow_average_mode_mean_condition} 
\end{figure}

\begin{figure}[ht]
\centering
\begin{minipage}[c]{0.48\textwidth}
\centering
    \includegraphics[width=\linewidth]{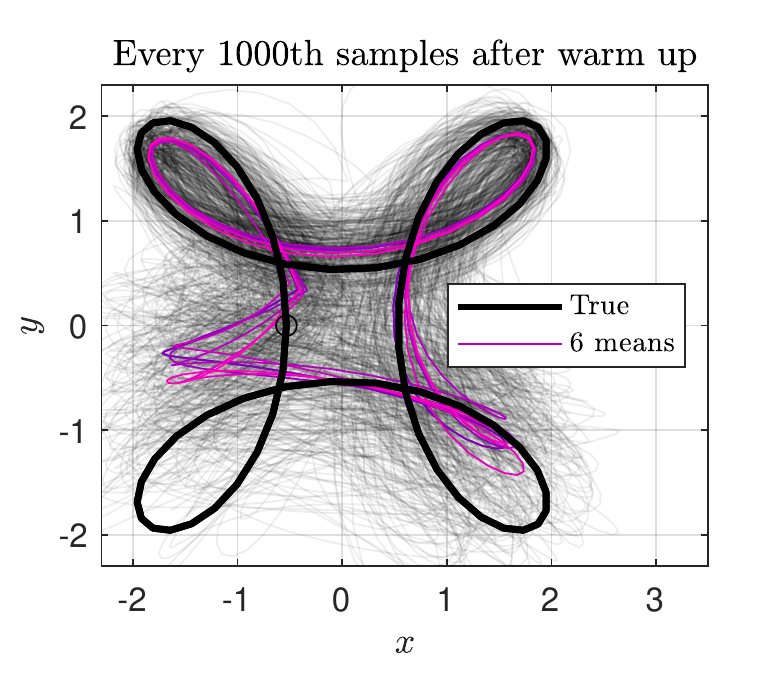}
    \end{minipage}
    \begin{minipage}[c]{0.48\textwidth}
\centering
        \includegraphics[width=\linewidth]{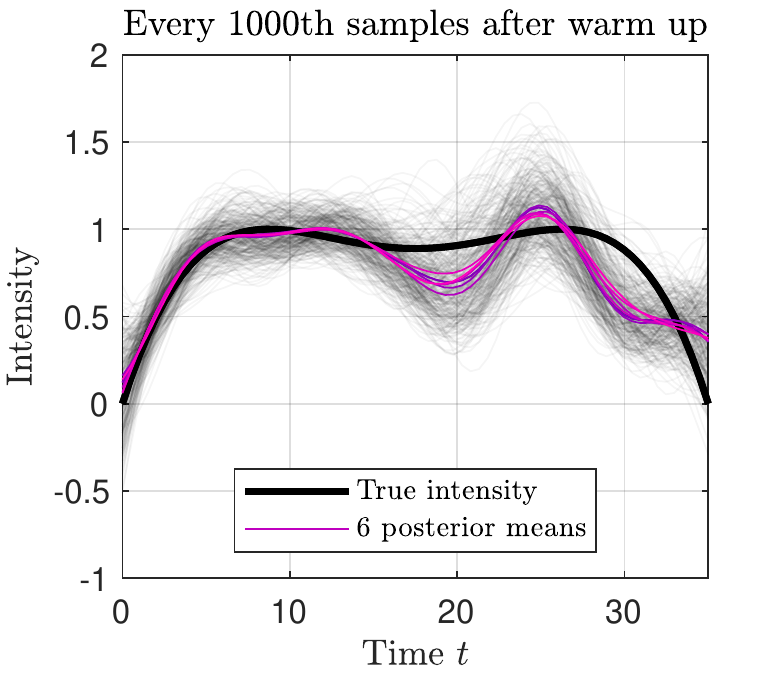}
\end{minipage}
\caption{Every 1000th sample from all 6 chains, and their posterior mean.}
\label{fig:bow_every_1000th_condition_half} 
\begin{minipage}[c]{0.45\textwidth}
\centering
    \includegraphics[width=\linewidth]{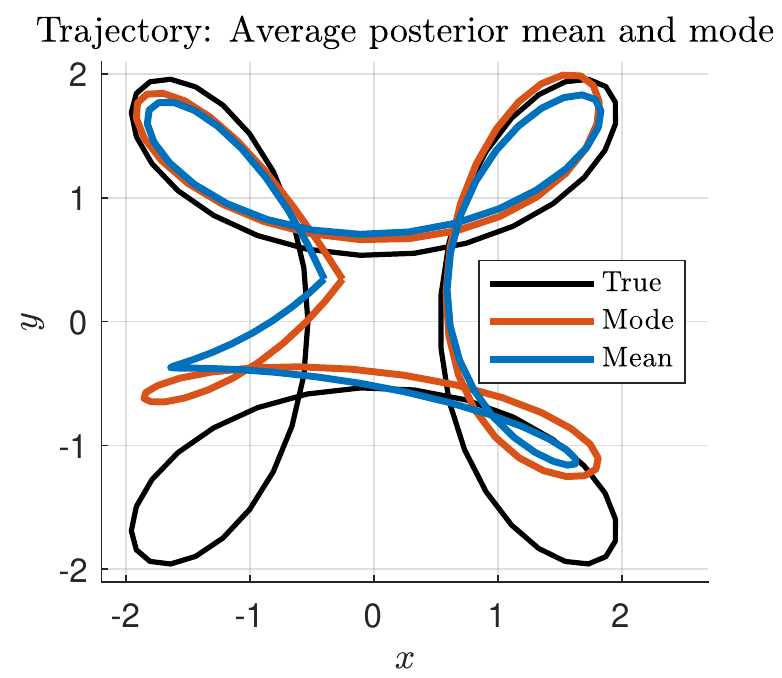}
\end{minipage}
\begin{minipage}[c]{0.45\textwidth}
\centering
    \includegraphics[width=\linewidth]{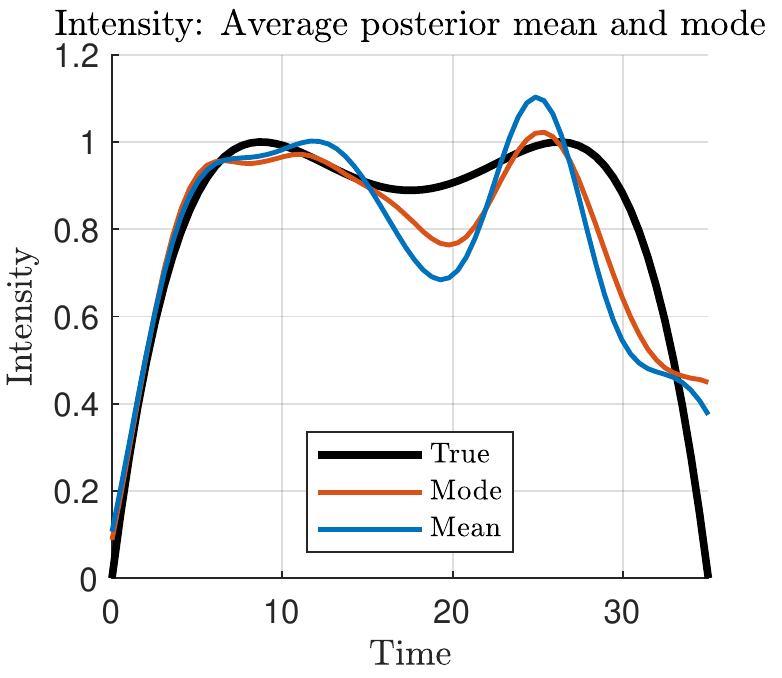}
\end{minipage}
\caption{The average posterior mean and mode of the trajectory and intensity.}
\label{fig:bow_average_mode_mean_condition_half} 
\end{figure}

We now consider a more complicated and longer trajectory with multiple crossing points given in parametric form as
\[
p_x(t)=\frac{8}{5}\sin(2\pi t/T_0+3\pi/2)+\frac{16}{15}\cos(6\pi t/T_0),\quad t\in[0,T_0],
\]
\[
p_y(t)=\frac{8}{5}\cos(2\pi t/T_0+3\pi/2)+\frac{16}{15}\sin(6\pi t/T_0),\quad t\in[0,T_0],
\]
\[
p_z(t)=0,\quad t\in[0,T]
\]
with the same emission intensity considered earlier
\begin{equation}
    q(t)=\begin{cases}
    -28.44(t/T_0)^4+56.89(t/T_0)^3-39.11(t/T_0)^2+10.67(t/T_0), &t\in[0,T_0] \\
    0, & \text{otherwise}
    \end{cases}.
\end{equation}
We fix the hyperparameters $\kappa=1.2$,  $\ell=4$. The trajectory and intensity with a number of initial draws from the GP prior are shown in Fig. \ref{fig:bow_set_up}.
%
Figure \ref{fig:bow_every_1000th} shows every 1000th sample from the pCN-MCMC chain together with the posterior means. We obtain a good reconstruction of the intensity while the trajectory deviates near the trajectory endpoints due to the low intensity near $t = 0$ and $t = T_0$. A comparison of the posterior mean and mode averaged over 6 chains is shown in Fig. \ref{fig:bow_average_mode_mean} and we see that while the averaged mode is slightly better than the averaged mean, the endpoint problem persists. This is reflected in the average trajectory error in Tbl. \ref{tab:bow_mean_error} compared to the previous cases.
\begin{table}[ht]
\begin{center}
  \begin{tabular}{ l | c | c |  }
    & Average mean & Average Mode \\ \hline
wavefield error & $2.0\cdot 10^{-3} $ & $1.7\cdot 10^{-3} $ \\ \hline
Average trajectory error & $2.3\cdot 10^{-1} $ & $2.2\cdot 10^{-1} $  \\ \hline
Average intensity error & $1.7 \cdot 10^{-2} $ & $2.0 \cdot 10^{-2} $  \\ \hline
  \end{tabular}
  \caption{Error of the average posterior mean and mode.}
  \label{tab:bow_mean_error}
\end{center}
\end{table}
We now show how the reconstruction can be greatly aided by incoorporating the prior knowledge that the trajectory is closed as discussed in Sec. \ref{sec:prior_knowledge}. Note that we use no knowledge of the exact trajectory coordinates at any time, we only condition the prior on $p(0) = p(T_0)$. This leads to the reconstruction in Fig. \ref{fig:bow_average_mode_mean_condition}.
%
%
%
%
%
Visually, the reconstruction of the trajectory has significantly improved: The prior provides information where the reconstruction previously suffered due to the low intensity. We still see some deviation near the endpoints and reconstruction could be further improved by conditioning on a smooth joining of the endpoint, i.e. $p'(0) = p'(T_0)$. Tbl. \ref{tab:bow_mean_error_condition_or_not} shows how the trajectory reconstruction has improved through the simple single-point conditioning. The ability to incoorporate prior information in such a flexible way is a key advantage of this method compared to e.g. algebraic methods which rely on point-wise estimation.
\begin{table}[ht]
\begin{center}
  \begin{tabular}{ l | c | c |  }
    & Not conditioned & Conditioned \\ \hline
wavefield error & $2.0\cdot 10^{-3} $ & $1.7\cdot 10^{-3} $ \\ \hline
Average trajectory error & $23\cdot 10^{-2} $ & $11\cdot 10^{-2} $  \\ \hline
Average intensity error & $1.7 \cdot 10^{-2} $ & $2.5 \cdot 10^{-2} $  \\ \hline
  \end{tabular}
  \caption{Errors using the average posterior mean from conditioned and unconditioned model.}
  \label{tab:bow_mean_error_condition_or_not}
\end{center}
\end{table}
%
Lastly, we show how the choice of measurement surface influences the reconstruction. We again condition on a closed trajectory but remove half the measurement points, leaving 213 sensors in a quarter sphere. The resulting reconstruction is shown in Fig. \ref{fig:bow_average_mode_mean_condition_half}.
The impact on the reconstruction is evident: The reconstruction suffers significantly in the lower half while we retain a good reconstruction in the upper half. A similar effect is observed in the intensity reconstruction, where the time interval $[10, 20]$ corresponds to the particle travelling in the lower half of the trajectory.

\FloatBarrier
\subsection{Case 4}

\begin{figure}[ht]
\centering
\begin{minipage}[c]{0.45\textwidth}
\centering
\includegraphics[width=\linewidth]{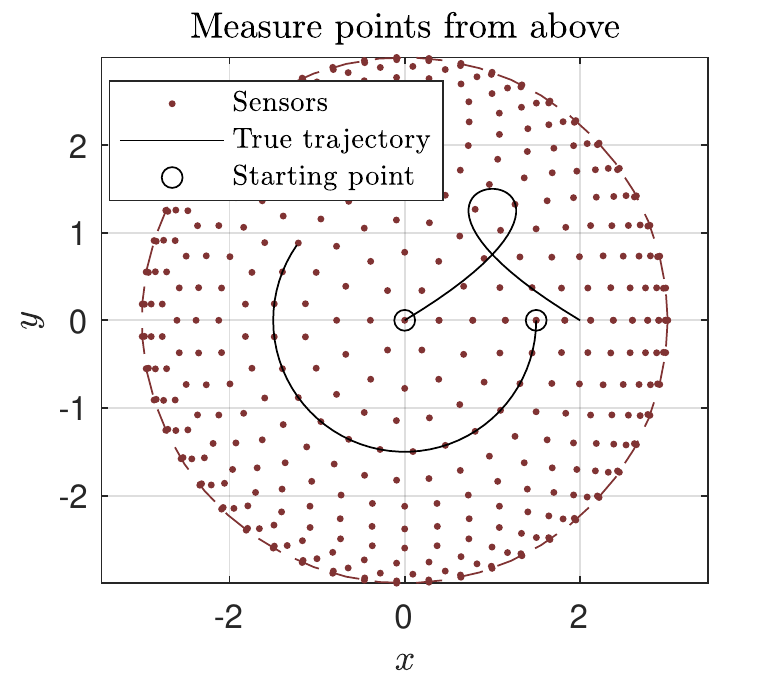}
$(a)$
\end{minipage}
\begin{minipage}[c]{0.45\textwidth}
\centering
\includegraphics[width=\linewidth]{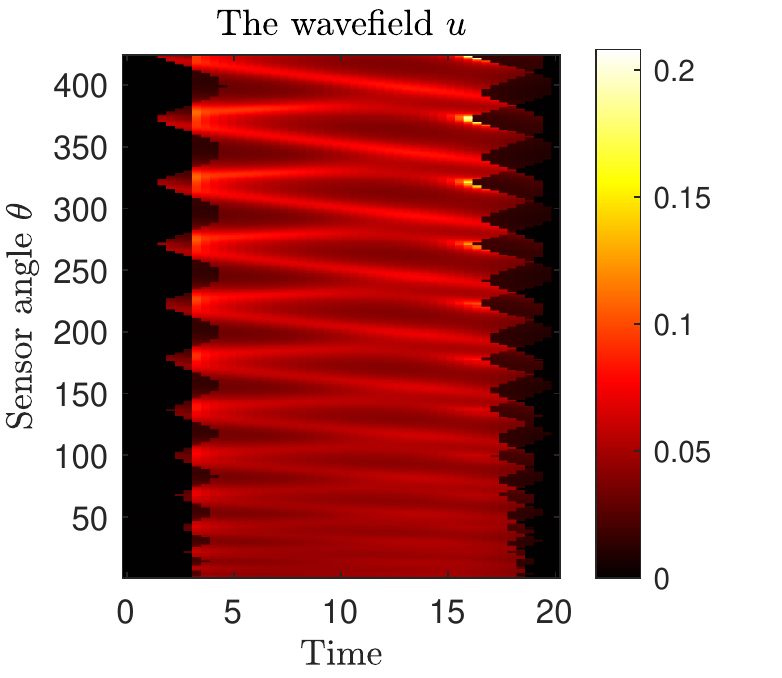}
$(b)$
\end{minipage}
\caption{The setup: $(a)$ shows the trajectories and sensor placement from above and $(b)$ shows the measurement data. }
\label{fig:two_sources_setup}
\end{figure}

\begin{figure}[ht]
\centering
\begin{minipage}[c]{0.45\textwidth}
\centering
\includegraphics[width=\linewidth]{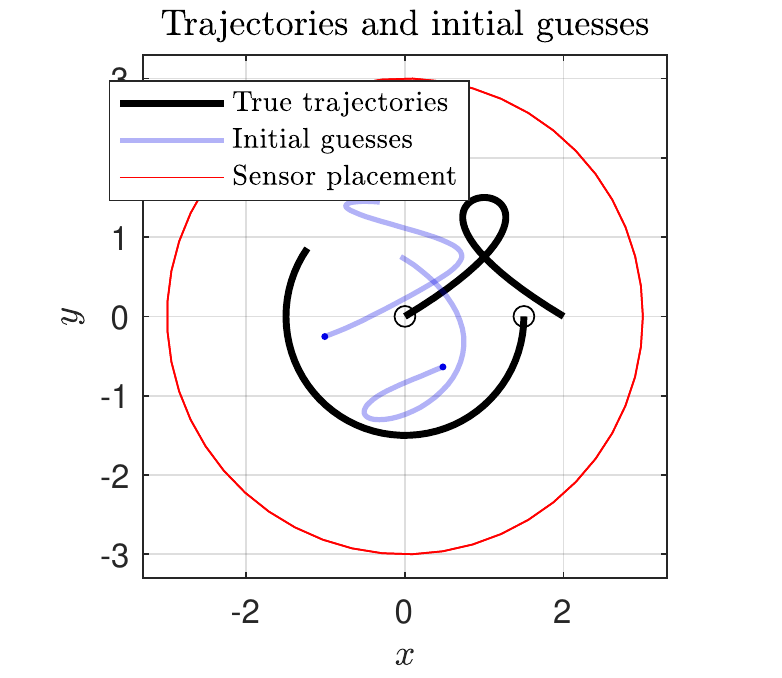}
$(a)$
\end{minipage}
\begin{minipage}[c]{0.45\textwidth}
\centering
\includegraphics[width=\linewidth]{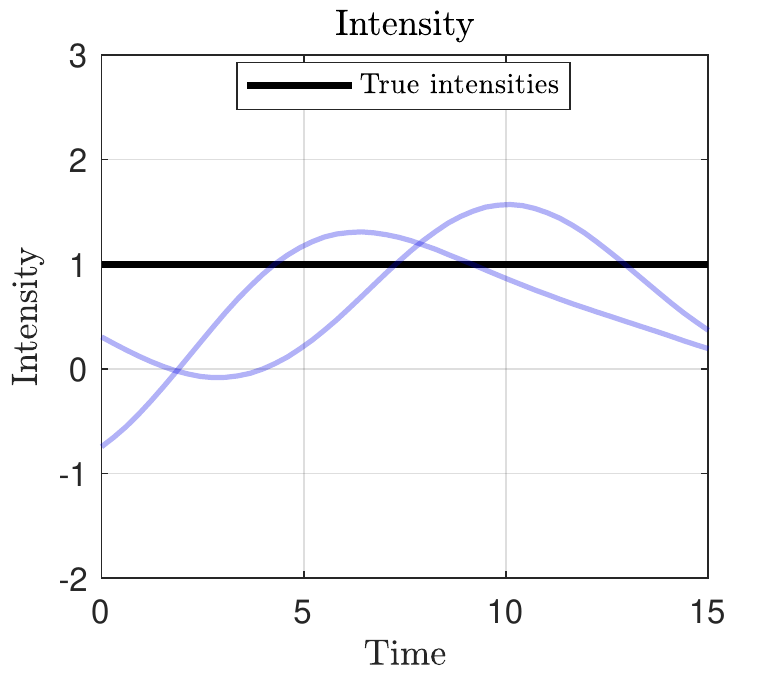}
$(b)$
\end{minipage}
\caption{Initial guesses of $(a)$ the trajectories and $(b)$ the intensity for two sources.}
\label{fig:two_sources_initial}
\end{figure}

\begin{figure}[ht]
\centering
\begin{minipage}[c]{0.45\textwidth}
\centering
\includegraphics[width=\linewidth]{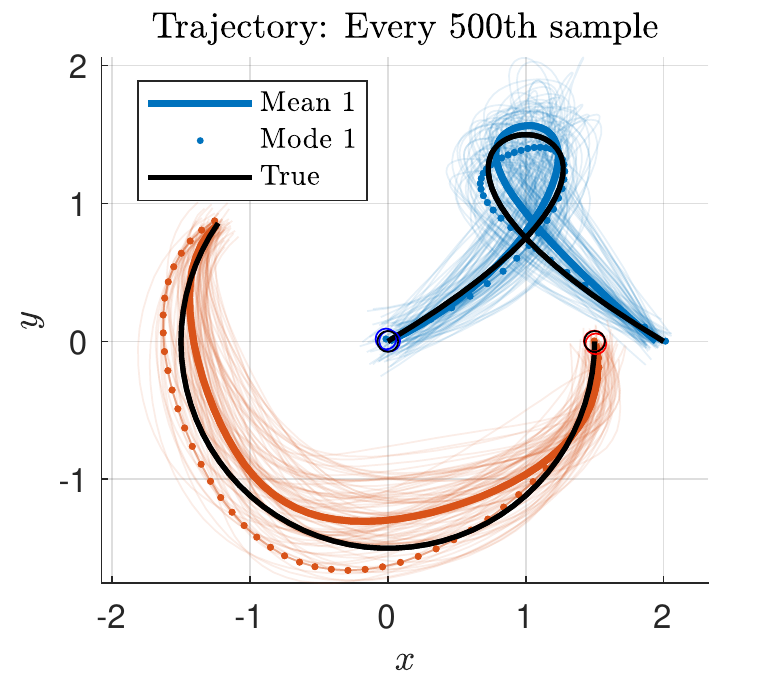}
$(a)$
\end{minipage}
\begin{minipage}[c]{0.45\textwidth}
\centering
\includegraphics[width=\linewidth]{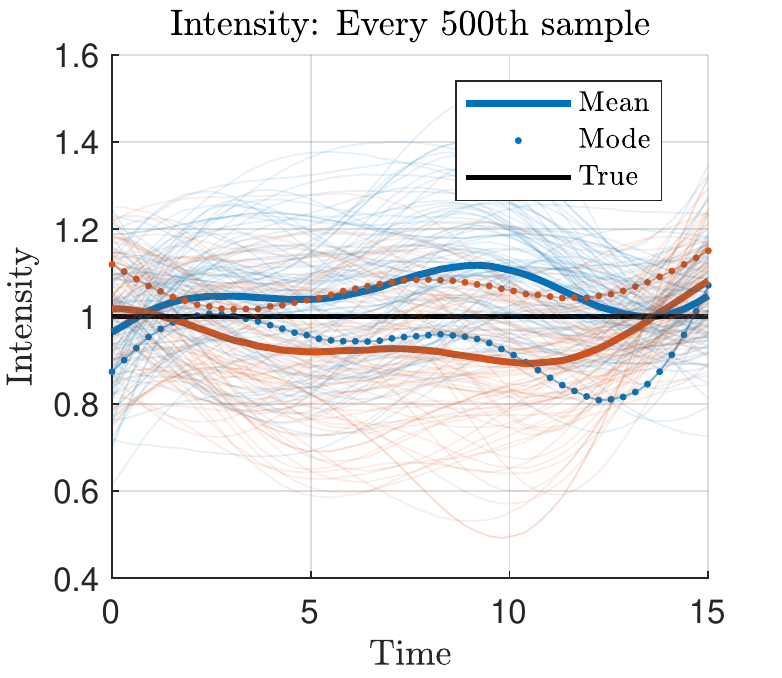}
$(b)$
\end{minipage}
\caption{Posterior mean, mode and every 500th sample of $(a)$ the trajectories and $(b)$ the intensities. The coloured solid lines are the posterior means and the dotted lines are the posterior modes.}
\label{fig:two_sources_samples}
\end{figure}

As a final example we demonstrate how our method performs in the presence of two simultaneously radiating sources. We consider the trajectories given in parametric form by
\begin{subequations}
\begin{align}
p_1(t)&=\left(16(t/T_0)^3-24(t/T_0)^2+5t/T_0,3(t/T_0)(1-t/T_0),0\right),& t&\in[0,T_0] \\
p_2(t)&=\left(1.5\cos(0.25t),-1.5\sin(0.25t),0\right), & t&\in[0,T_0]
\end{align}
\end{subequations}
with constant emission intensity
\[
q_1(t)=q_2(t)=1,\,\,\,t\in[0,T_0],\quad q_1(t)=q_2(t)=0,\,\,\,t>T_0
\]
and we choose the hyperparameters $\kappa = 1$ and $\ell = 5$ for all priors for the reconstruction. Figure \ref{fig:two_sources_setup} shows the trajectories and the distribution of 424 sensors over a hemisphere together with the measured wave field. Figure \ref{fig:two_sources_initial} shows the trajectories and intensities together with two draws from the prior distribution.
We then generate 100.000 samples from the MCMC algorithm and the resulting Bayesian reconstruction is shown in Fig. \ref{fig:two_sources_samples}, showing the posterior mean and mode together with every 500th sample from the MCMC chain. While the reconstruction still yields reasonable results, we see that the reconstruction suffers slightly when another emitter is introduced compared to the single-emitter case. 
\FloatBarrier 

\bibliographystyle{plain}
\bibliography{bib}

\begin{thebibliography}{10}

\bibitem{antil2021novel}
Harbir Antil, Howard~C Elman, Akwum Onwunta, and Deepanshu Verma.
\newblock Novel deep neural networks for solving bayesian statistical inverse.
\newblock {\em arXiv preprint arXiv:2102.03974}, 2021.

\bibitem{baorecovering}
Gang Bao, Yuantong Liu, and Faouzi Triki.
\newblock Recovering point sources for the inhomogeneous helmholtz equation.
\newblock {\em Inverse Problems}, 37(9):095005, 2021.

\bibitem{bao2021recovering}
Gang Bao, Yuantong Liu, and Faouzi Triki.
\newblock Recovering simultaneously a potential and a point source from cauchy
  data.
\newblock {\em Minimax Theory and its Applications}, 06(2):227--238, 2021.

\bibitem{beskos2017geometric}
Alexandros Beskos, Mark Girolami, Shiwei Lan, Patrick~E Farrell, and Andrew~M
  Stuart.
\newblock Geometric mcmc for infinite-dimensional inverse problems.
\newblock {\em Journal of Computational Physics}, 335:327--351, 2017.

\bibitem{bruckner2000determination}
Gottfried Bruckner and Masahiro Yamamoto.
\newblock Determination of point wave sources by pointwise observations:
  stability and reconstruction.
\newblock {\em Inverse problems}, 16(3):723, 2000.

\bibitem{conrad2016accelerating}
Patrick~R Conrad, Youssef~M Marzouk, Natesh~S Pillai, and Aaron Smith.
\newblock Accelerating asymptotically exact mcmc for computationally intensive
  models via local approximations.
\newblock {\em Journal of the American Statistical Association},
  111(516):1591--1607, 2016.

\bibitem{Duistermaat}
H.~Duistermaat.
\newblock {\em Distributions: Theory and Applications}.
\newblock Birkh\"auser, 2010.

\bibitem{el2001determination}
Abdellatif El~Badia and T~Ha-Duong.
\newblock Determination of point wave sources by boundary measurements.
\newblock {\em Inverse Problems}, 17(4):1127, 2001.

\bibitem{gelman2013bayesian}
Andrew Gelman, John~B Carlin, Hal~S Stern, and Donald~B Rubin.
\newblock {\em Bayesian data analysis}.
\newblock Chapman and Hall/CRC, 2013.

\bibitem{AlJebawy}
A.~Elbadia H.~Al~Jebawy and F.~Triki.
\newblock Inverse moving point source problem for the wave equation.
\newblock {\em Arxiv}, 2022.

\bibitem{HI}
L.~H\"ormander.
\newblock {\em The Analysis of Linear Partial Differential Operators I}.
\newblock 2003.

\bibitem{hu2}
G.~Hu, Y.~Kian, P.~Li, and Y.~Zhao.
\newblock Inverse moving source problems in electrodynamics.
\newblock {\em Inverse Problems}, 35:075001, 2019.

\bibitem{hu}
G.~Hu, Y.~Liu, and M.~Yamamoto.
\newblock Inverse moving source problem for fractional diffusion (-wave)
  equations: Determination of orbits.
\newblock In {\em International Conference on Inverse Problems}, pages 81--100.
  Springer, 2018.

\bibitem{Jackson}
J.~D. Jackson.
\newblock {\em Classical electrodynamics}.
\newblock 1999.

\bibitem{komornik2002upper}
Vilmos Komornik and Masahiro Yamamoto.
\newblock Upper and lower estimates in determining point sources in a wave
  equation.
\newblock {\em Inverse Problems}, 18(2):319, 2002.

\bibitem{komornik2005estimation}
Vilmos Komornik and Masahiro Yamamoto.
\newblock Estimation of point sources and applications to inverse problems.
\newblock {\em Inverse Problems}, 21(6):2051, 2005.

\bibitem{nakaguchi}
E.~Nakaguchi, H.~Inui, and K.~Ohnaka.
\newblock An algebraic reconstruction of a moving point source for a scalar
  wave equation.
\newblock {\em Inverse Problems}, 28:065018, 2012.

\bibitem{nakaguchi2012algebraic}
Etsushi Nakaguchi, Hirokazu Inui, and Kohzaburo Ohnaka.
\newblock An algebraic reconstruction of a moving point source for a scalar
  wave equation.
\newblock {\em Inverse Problems}, 28(6):065018, 2012.

\bibitem{Nedelec}
J.-C. N\'ed\'elec.
\newblock {\em Acoustic and Electromagnetic Equations: Integral Representations
  for Harmonic Problems}.
\newblock Springer, 2001.

\bibitem{ohe}
T.~Ohe.
\newblock Real-time reconstruction of moving point/dipole wave sources from
  boundary measurements.
\newblock {\em Inverse Problems in Science and Engineering}, 28:1057--1102,
  2020.

\bibitem{ohe2011real}
Takashi Ohe, Hirokazu Inui, and Kohzaburo Ohnaka.
\newblock Real-time reconstruction of time-varying point sources in a
  three-dimensional scalar wave equation.
\newblock {\em Inverse Problems}, 27(11):115011, 2011.

\bibitem{rashedi2015stable}
Kamal Rashedi and Mourad Sini.
\newblock Stable recovery of the time-dependent source term from one
  measurement for the wave equation.
\newblock {\em Inverse Problems}, 31(10):105011, 2015.

\end{thebibliography}


\end{document}